\documentclass[3p,times]{elsarticle}

\usepackage{amsfonts}
\usepackage{mathrsfs}
\usepackage{amsmath}
\usepackage{amssymb}

 \usepackage{lineno}
\newtheorem{theorem}{\indent Theorem}[section]
\newtheorem{corollary}{\indent Corollary}[section]

\newtheorem{definition}{\indent Definition}[section]
\newtheorem{lemma}{\indent Lemma}[section]
\newtheorem{remark}{\indent Remark}[section]

\journal{Communications in Mathematical Sciences}

\begin{document}

\begin{frontmatter}



\title{Finite dimensional global attractor of the Cahn-Hilliard-Navier-Stokes system with dynamic boundary conditions}


\author{Bo You$^{a,*},$ Fang Li$^b,$ Chang Zhang$^c$}

\address{$^a$ School of Mathematics and Statistics, Xi'an Jiaotong University, Xi'an, 710049, P. R. China\\
$^b$ School of Mathematics and Statistics, Xidian University, Xi'an, 710126, P. R. China\\
$^c$ Department of Mathematics, Nanjing University, Nanjing, 210093, P. R. China
}

 \cortext[cor1]{Corresponding author.\\
 E-mail addresses:youb2013@xjtu.edu.cn(B. You),fli@xidian.edu.cn (F. Li), chzhnju@126.com (C. Zhang)}
\begin{abstract}
In this paper, we mainly consider the long-time behavior of solutions for the Cahn-Hilliard-Navier-Stokes system with dynamic boundary conditions and two polynomial growth nonlinearities of arbitrary order. We prove the existence of a finite dimensional global attractor for the Cahn-Hilliard-Navier-Stokes system with dynamic boundary conditions by using the $\ell$-trajectories method.
\end{abstract}


\end{frontmatter}

\section{Introduction}
\def\theequation{1.\arabic{equation}}\makeatother
\setcounter{equation}{0}
In this paper, we consider the following Cahn-Hilliard-Navier-Stokes system:
\begin{equation}\label{1.1}
\begin{cases}
&\frac{\partial u}{\partial t}-\nu \Delta u+(u\cdot\nabla)u+\nabla p+\lambda\phi\nabla\mu=h(x),\,\,\,(x,t)\in\Omega\times\mathbb{R}^+,\\
&\nabla\cdot u=0,\,\,\,(x,t)\in\Omega\times\mathbb{R}^+,\\
&\frac{\partial \phi}{\partial t}+u\cdot\nabla\phi-\gamma\Delta\mu=0,\,\,\,(x,t)\in\Omega\times\mathbb{R}^+,\\
&\mu=-\Delta\phi+f(\phi),\,\,\,(x,t)\in\Omega\times\mathbb{R}^+.
\end{cases}
\end{equation}
Equation \eqref{1.1} is subject to the following dynamic boundary conditions
\begin{equation}\label{1.2}
\begin{cases}
&u(x,t)=0,\,\,\,(x,t)\in\Gamma\times\mathbb{R}^+,\\
&\frac{\partial \mu}{\partial \vec{n}}=0,\,\,\,(x,t)\in\Gamma\times\mathbb{R}^+,\\
&\frac{\partial \phi}{\partial t}=\alpha\Delta_\Gamma\phi-\frac{\partial \phi}{\partial \vec{n}}-\beta\phi-g(\phi),\,\,\,(x,t)\in\Gamma\times\mathbb{R}^+
\end{cases}
\end{equation}
 and initial conditions
\begin{equation}\label{1.3}
\begin{cases}
&u(x,0)=u_0(x),\,\,\,x\in\Omega,\\
&\phi(x,0)=\phi_0(x),\,\,\,x\in\Omega,\\
&\phi(x,0)=\theta_0(x),\,\,\,x\in\Gamma,
\end{cases}
\end{equation}
where $\Omega\subset\mathbb{R}^2$ is a bounded domain with smooth boundary $\Gamma$ and $\mathbb{R}^+=[0,+\infty),$ $\nu>0$ is the viscosity, $\lambda>0$ is a surface tension parameter, $\alpha>0,$ $\beta>0$ are constants, $\gamma>0$ is the elastic relaxation time, $h(x)=(h_1(x),h_2(x))$ is the external force, $u(x,t)=(u_1(x,t),u_2(x,t))$ denotes the average velocity and $\phi$ is the difference of the two fluid concentrations, $p$ is the fluid pressure, $\vec{n}$ is the unit external normal vector on $\Gamma,$ $\Delta_\Gamma$ is the Laplace-Beltrami operator on the surface $\Gamma$ of $\Omega.$

 To study problem \eqref{1.1}-\eqref{1.3}, we assume the following conditions:

$(H_1)$ the function $f\in C^1(\mathbb{R},\mathbb{R})$ satisfies that there exists a positive constant $C_1$ such that
\begin{align}\label{1.4}
|f'(u)-f'(v)|\leq C_1|u-v|(|u|^{p-3}+|v|^{p-3}+1)
\end{align}
for any $u,v\in\mathbb{R}$ and
\begin{align}\label{1.5}
 &c_1|u|^p-k_1\leq f(u)u\leq c_2|u|^p+k_1,
\end{align}
where $c_i> 0$ ($i=1,2$), $p\geq 3,$ $ k_1> 0.$

$(H_2)$ the function $g\in C(\mathbb{R},\mathbb{R})$ satisfies that there exists a positive constant $C_2$ such that
\begin{align}\label{1.6}
|g(u)-g(v)|\leq C_2|u-v|(|u|^{q-2}+|v|^{q-2}+1)
\end{align}
for any $u,v\in\mathbb{R}$ and
\begin{align}\label{1.7}
 &c_3|u|^q-k_2\leq g(u)u\leq c_4|u|^q+k_2,
\end{align}
where $c_i> 0$ ($i=3,4$), $q>2,$ $k_2>0.$

 Dynamic boundary conditions were recently proposed by physicists to describe spinodal decomposition of binary mixtures where the effective interaction between the wall (i.e., the boundary) and two mixture components is short-ranged, and this type of boundary conditions is very natural in many mathematical
 models such as heat transfer in a solid in contact with a moving
 fluid, thermoelasticity, diffusion phenomena, heat transfer in two
 medium, problems in fluid dynamics. The well-posedness and long-time behavior of solutions for many equations with dynamical boundary conditions have been studied extensively(see \cite{ajm, ajm1, bar, cr, ca, ca1, fzh, gcg5, gcg6, gcg7, ma, pj, pj1, pl, rr, yl}). For example, the global well-posedness of solutions for the non-isothermal Cahn-Hilliard equation with dynamic boundary conditions was proved in \cite{gcg4}. In \cite{gcg5}, the author proved the existence and uniqueness of a global solution for a Cahn-Hilliard model in bounded domains with permeable walls. The global existence and uniqueness of solutions for the Cahn-Hilliard equation with highest-order boundary conditions were proved in \cite{rr}. In \cite{pj1}, the authors proved the maximal regularity and asymptotic behavior of solutions for the Cahn-Hilliard equation with dynamic boundary conditions. The fact that any global weak/strong solution of the Cahn-Hilliard equation with dynamic boundary conditions converges to a single steady state as time $t\rightarrow +\infty$ was proved in \cite{cr}. In \cite{gcg8}, the author proved the existence of a global attractor and an exponential attractor in $H^1(\Omega)$ for a homogeneous two-phase flow model and established any global weak/strong solution converges to a single steady state as time $t\rightarrow +\infty,$ and provided its convergence rate. In \cite{gcg6}, the author proved the existence of an exponential attractor for a Cahn-Hilliard model in bounded domains with permeable walls. The existence of a global attractor for the reaction-diffusion equation with dynamical boundary conditions was proved in \cite{fzh}. In \cite{ma}, the authors proved the existence of an exponential attractor for the Cahn-Hilliard equation with dynamical boundary conditions. In \cite{yb}, the authors proved the existence of a global attractor for $p$-Laplacian equations with dynamical boundary conditions by using asymptotical a priori estimates.The well-posedness of solutions and the existence of a global attractor of the Cahn-Hilliard-Brinkman system with dynamic boundary conditions was proved in \cite{yb1}.

Diffuse-interface methods in fluid mechanics are widely used by many researchers in order to describe the behavior of complex fluids (see \cite{adm}).
A diffuse interface variant of Cahn-Hilliard-Navier-Stokes system has been proposed to
model the motion of an isothermal mixture of two immiscible and incompressible fluids subject to
phase separation (see\cite{gme, hpc, jd}). The coupled system consists of a convective Cahn-Hilliard equation for the order parameter, i.e., the difference of
the relative concentrations of the two phases, coupled with the Navier-Stokes
equations for the average fluid velocity. The Cahn-Hilliard-Navier-Stokes system has been
investigated from the numerical (see \cite{fxb, kd, kd1}) and analytical (see, e.g., \cite{ah1, ah,bs, bf, ccs, gcg, gcg1, gcg2, gcg9, svn,tm, zl, zy}) viewpoint in several papers. The long-time behavior and well-posedness of solutions for the two dimensional Cahn-Hilliard-Navier-Stokes system were proved in \cite{gcg}. Thanks to the shortage of the uniqueness of solutions, the authors have proved the existence of trajectory attractors for binary fluid mixtures in 3D in \cite{gcg1}. In \cite{gcg2}, the authors have considered the instability of two-phase flows and provided a lower bound on the dimension of the global attractor of the Cahn-Hilliard-Navier-Stokes system. The existence of pullback exponential attractor for a two dimensional Cahn-Hilliard-Navier-Stokes system in \cite{bs}. In \cite{tm}, the author has proved the existence of pullback attracots for a two dimensional non-autonomous Cahn-Hilliard-Navier-Stokes system. Recently, the authors have considered the Cahn-Hilliard-Navier-Stokes system with moving contact lines and proved any suitable global energy solution will convergent to a single equilibrium in \cite{gcg9}. In \cite{cl}, the authors have proved the well-posedness of solutions for the viscous Cahn-Hilliard-Navier-Stokes system with dynamic boundary conditions and considered the regularity of the weak solutions under some additional assumptions that $\phi_t(0)\in H^1(\bar{\Omega},d\sigma)$ and $u_t(0)\in H.$ However, to the best of our knowledge, there are no results related to the existence of a finite dimensional  global attractor for the dissipative dynamical system with dynamical boundary conditions.

In this paper, we will consider the well-posedness and the long-time behavior of solutions for the Cahn-Hilliard-Navier-Stokes system with dynamical boundary conditions and a polynomial growth nonlinearity of arbitrary order. When we consider the long-time behavior of solutions for the Cahn-Hilliard-Navier-Stokes system with dynamic boundary conditions, there are two difficulties: first of all, comparing to Cahn-Hilliard equation with dynamic boundary conditions, since the the coupled term arises and the additional assumptions that $\phi_t(0)\in H^1(\bar{\Omega},d\sigma)$ and $u_t(0)\in H$ specified in Lemma 2.3 of \cite{cl} can not be obtained for the weak solution at sufficiently large time in general such that we cannot obtain the existence of an absorbing set for problem \eqref{1.1}-\eqref{1.7} in a more regular phase space than $H\times V_I.$ Secondly, comparing to the Cahn-Hilliard-Navier-Stokes system with Neumann boundary conditions, thanks to
\begin{align*}
\int_\Omega\phi_t\Delta^2\phi\,dx=\int_\Gamma\phi_t\frac{\partial\Delta\phi}{\partial\vec{n}}-\int_\Gamma\Delta\phi\frac{\partial\phi_t}{\partial\vec{n}}+\cdots,
\end{align*}
it is very tricky to deal with these two terms on the right hand side such that we cannot choose $\Delta^2\phi$ as a test function to prove the smooth property of the difference of two solutions and the differentiability of the corresponding semigroup on the global attractor for problem \eqref{1.1}-\eqref{1.7}. Therefore, the standard scheme of estimating the fractal dimension of the global attractor does not work. To overcome this difficulty, inspired by the idea of the method of $\ell$-trajectories for any small $\ell>0$ proposed in \cite{mj}, in this paper, we first define a semigroup $\{L_t\}_{t\geq 0}$ on some subset $X_\ell$ of $L^2(0,\ell;H\times V_I)$ induced by the semigroup $\{S_I(t)\}_{t\geq0}$ generated by problem \eqref{1.1}-\eqref{1.7}, and then, we prove the existence of a global attractor $\mathcal{A}_\ell$ in $X_\ell$ for the semigroup $\{L_t\}_{t\geq 0}$ by the method of $\ell$-trajectories and estimate the fractal dimension of the global attractor by using the smooth property of the semigroup $\{L_t\}_{t\geq 0}.$ Finally, by defining a Lipschitz continuous operator on the global attractor $\mathcal{A}_\ell,$ we obtain the existence of a finite dimensional global attractor $\mathcal{A}$ in the original phase space $H\times V_I$ for problem \eqref{1.1}-\eqref{1.7}.

Throughout this paper, let $C$ be a generic constant that is independent of the initial datum of $(u,\phi).$ Define the average of function $\phi(x)$ over $\Omega$ as
\begin{align*}
m\phi=\frac{1}{|\Omega|}\int_\Omega \phi(x)\,dx.
\end{align*}
\section{\bf Preliminaries}
\def\theequation{2.\arabic{equation}}\makeatother
\setcounter{equation}{0}
In order to study the problem \eqref{1.1}-\eqref{1.7}, we introduce the space of divergence-free functions defined by
\begin{align*}
\mathcal{V}=\{u\in (\mathcal{C}_c^{\infty}(\Omega))^2:\nabla\cdot u=0\}.
\end{align*}
Denote by $H$ and $V$ the closure of $\mathcal{V}$ with respect to the norms in $(L^2(\Omega))^2$ and $(H_0^1(\Omega))^2,$ respectively.

We define the Lebesgue spaces as follows
\begin{align*}
L^p(\Gamma)=\left\{v:\|v\|_{L^p(\Gamma)}<\infty\right\},
\end{align*}
where
\begin{align*}
\|v\|_{L^p(\Gamma)}=\left(\int_{\Gamma}|v|^p\,dS\right)^{\frac{1}{p}}
\end{align*}
for $p\in [1,\infty).$ Moreover, we have
\begin{align*}
L^p(\Omega)\oplus L^q(\Gamma)=L^{p,q}(\bar{\Omega},d\sigma),\,\,\,\,\,p,q\in [1,\infty)
\end{align*}
and
\begin{align*}
\|U\|_{L^{p,q}(\bar{\Omega},d\sigma)}=(\int_{\Omega}|u|^p\,dx)^{\frac{1}{p}}+(\int_{\Gamma}|v|^q\,dS)^{\frac{1}{q}}
\end{align*}
for any $U=(u,v)\in L^{p,q}(\bar{\Omega},d\sigma),$ where the measure $d\sigma=dx|_{\Omega}\oplus dS|_{\Gamma}$ on $\bar{\Omega}$ is defined by $\sigma(A)=|A\cap\Omega|+S(A\cap \Gamma)$ for any measurable set $A\subset \bar{\Omega}.$

We also define the Sobolev space $H^1(\bar{\Omega},d\sigma)$ as the closure of $C^1(\overline{\Omega})$ with respect to the norm given by
\begin{align*}
\|\phi\|_{H^1(\bar{\Omega},d\sigma)}^2=\left(\int_{\Omega}|\nabla \phi|^2\,dx+\int_{\Gamma}\alpha|\nabla_\Gamma \phi|^2+\beta|\phi|^2\,dS\right)^{\frac{1}{2}}
\end{align*}
for any $\phi\in C^1(\overline{\Omega}),$  denote by $X^*$ the dual space of $X$ and let $H^s(\Omega),$ $H^s(\Gamma)$ $(s\in\mathbb{R})$ be the usual Sobolev spaces. In general, any vector $\theta\in L^p(\bar{\Omega},d\sigma)$ will be of the form $(\theta_1,\theta_2)$ with $\theta_1\in L^p(\Omega,dx)$ and $\theta_2\in L^p(\Gamma,dS),$ and there need not be any connection between $\theta_1$ and $\theta_2.$

Let the operator $A: H^1(\bar{\Omega},d\sigma)\rightarrow (H^1(\bar{\Omega},d\sigma))^{*} $ be associated with the bilinear form defined by
\begin{align}\label{2.1}
\langle A\phi, \psi\rangle
=\int_{\Omega}\nabla \phi\cdot\nabla\psi\,dx+\int_{\Gamma}\alpha\nabla_\Gamma \phi\cdot\nabla_\Gamma \psi+\beta\phi\psi\,dS
\end{align}
for any $\phi,$ $\psi\in H^1(\bar{\Omega},d\sigma).$

\begin{remark}(\cite{gcg3})\label{2.2}
$C(\bar{\Omega})$ is a dense subspace of $L^2(\bar{\Omega},d\sigma)$ and a closed subspace of $L^{\infty}(\bar{\Omega},d\sigma).$
\end{remark}

Next, we recall briefly some lemmas used to prove the well-posedness of weak solutions and the existence of a finite dimensional global attractor for problem \eqref{1.1}-\eqref{1.7}.
\begin{lemma}(\cite{rjc})\label{2.3}
{\rm Let $\mathcal {O}$ be a bounded domain in $\mathbb{R}^n$ and let $1<q<\infty.$
Assume that $\{g_{n}\}\subset L^{q}(\mathcal {O})$ with $\|\{g_n\}\|_{L^q(\mathcal
{O})}\leq C$, where $C$ is independent of $n$ and there exists $g\in L^q(\mathcal {O})$ such that $\{g_n\}\rightarrow
g$, as $n\rightarrow\infty$, almost everywhere in $\mathcal {O}.$ Then $ g_n\rightarrow g$, as
$n\rightarrow\infty$ weakly in $L^q(\mathcal {O})$.}
\end{lemma}

\begin{lemma}(\cite{gcg4, ma})\label{2.4}
Let $\Omega\subset\mathbb{R}^2$ be a bounded domain with smooth boundary $\Gamma.$ Consider the following linear problem
\begin{equation*}
\begin{cases}
&-\Delta \phi=j_1,\,\,\,x\in\Omega,\\
&-\alpha\Delta_\Gamma\phi+\frac{\partial\phi}{\partial\vec{n}}+\beta\phi=j_2,\,\,\,x\in\Gamma.
\end{cases}
\end{equation*}
Assume that $(j_1,j_2)\in H^s(\bar{\Omega},d\sigma),$ $s\geq 0,$ $s+\frac{1}{2}\not\in\mathbb{N}.$ Then the following estimate holds
\begin{align*}
\|\phi\|_{H^{s+2}(\bar{\Omega},d\sigma)}\leq C(\|j_1\|_{H^s(\Omega)}+\|j_2\|_{H^s(\Gamma)})
\end{align*}
for some constant $C>0.$
\end{lemma}

\begin{lemma}(\cite{tr})\label{2.5}
 Let $V,$ $H,$ $V^*$ be three Hilbert spaces such that $V\subset H=H^*\subset V^*,$ where $H^*$ and $V^*$ are the dual spaces of $H$ and $V,$ respectively. Suppose $u\in L^2(0,T;V)$ and $\frac{\partial u}{\partial t}\in L^2(0,T;V^*).$ Then $u$ is almost everywhere equal to a function continuous from $[0,T]$ into $H.$
\end{lemma}

\begin{lemma}(\cite{cvv,jn, mj, mj1,sj2})\label{2.6}
Assume that $p_1\in(1,\infty],$ $p_2\in [1,\infty).$ Let $X$ be a Banach space and let $X_0,$ $X_1$ be separable and reflexive Banach spaces such that $X_0\subset\subset X\subset X_1.$ Then
\begin{align*}
Y=\{u\in L^{p_1}(0,\ell;X_0):u'\in L^{p_2}(0,\ell;X_1)\}\subset\subset L^{p_1}(0,\ell;X),
\end{align*}
where $\ell$ is a fixed positive constant.
\end{lemma}

 \begin{definition}(\cite{rjc,tr})
Let $\{S(t)\}_{t\geq0}$ be a semigroup on a Banach space $X.$ A set $\mathcal{A}\subset X$ is said to be a global attractor if the following conditions hold:
\begin{itemize}
\item [(i)] $\mathcal{A}$ is compact in $X.$
\item [(ii)] $\mathcal{A}$ is strictly invariant, i.e., $S(t)\mathcal{A}=\mathcal{A}$ for any $t\geq 0.$
\item [(iii)] For any bounded subset $B\subset X$ and for any neighborhood $\mathcal{O}=\mathcal{O}(\mathcal{A})$ of $\mathcal{A}$ in $X,$
there exists a time $\tau_0=\tau_0(B)$ such that $S(t)B\subset\mathcal{O}(\mathcal{A})$ for any $t\geq \tau_0.$
\end{itemize}
\end{definition}
\begin{lemma}(\cite{mj1})\label{2.7}
Let $X$ be a (subset of) Banach space and $(S(t),X)$ be a dynamical system. Assume that there exists a compact set $K\subset X$ which is uniformly absorbing and positively invariant with respect to $S(t).$ Let moreover $S(t)$ be continuous on $K.$ Then $(S(t),X)$ has a global attractor.
\end{lemma}
\begin{definition}(\cite{rjc, tr})
Let $H$ be a separable real Hilbert space. For any non-empty compact subset $K\subset H,$ the fractal dimension of $K$ is the number
\begin{align*}
d_f(K)=\limsup_{\epsilon\rightarrow 0^+}\frac{\log(N_\epsilon(K))}{\log(\frac{1}{\epsilon})},
\end{align*}
where $N_\epsilon(K)$ denotes the minimum number of open balls in $H$ with radii $\epsilon>0$ that are necessary to cover $K.$
\end{definition}
\begin{lemma}(\cite{mj1})\label{2.8}
Let $X,$ $Y$ be norm spaces such that $X\subset\subset Y$ and $\mathcal{A}\subset Y$ be bounded. Assume that there exists a mapping $L$ such that $L\mathcal{A}=\mathcal{A}$ and $L:Y\rightarrow X$ is Lipschitz continuous on $\mathcal{A}.$ Then $d_f(\mathcal{A})$ is finite.
\end{lemma}
\begin{lemma}(\cite{mj1})\label{2.9}
Let $X$ and $Y$ be two metric spaces and $f:X\rightarrow Y$ be $\alpha$-H\"{o}lder continuous on the subset $A\subset X.$ Then
\begin{align*}
d_F(f(A),Y)\leq\frac{1}{\alpha}d_F(A,X).
\end{align*}
In particular, the fractal dimension does not increase under a Lipschitz continuous mapping.
\end{lemma}
Finally, we give the definition of weak solutions for problem \eqref{1.1}-\eqref{1.7}.
\begin{definition}
Assume that $h\in L^2(\Omega)$ and $(H_1)$-$(H_2)$ hold. For any $(u_0,\Phi_0)=(u_0,\phi_0,\theta_0)\in H\times H^1(\bar{\Omega},d\sigma)$ and any fixed $T>0,$ a function $(u,\phi)$ is called a weak solution of problem \eqref{1.1}-\eqref{1.7} on $(0,T),$ if
\begin{align*}
&\mu\in L^2(0,T;H^1(\Omega))\,\,\,\text{is given by the fourth equation of \eqref{1.1}}
\end{align*}
and
\begin{align*}
&\phi\in\mathcal{C}([0,T];H^1(\bar{\Omega},d\sigma))\cap L^2(0,T;H^2(\bar{\Omega},d\sigma)),\\
&u\in \mathcal{C}([0,T];H)\cap L^2(0,T;V),\\
&(u_t,\phi_t)\in L^2(0,T;V^*\times (H^1(\bar{\Omega},d\sigma))^*)
\end{align*}
satisfy
\begin{align*}
&\int_{\Omega}u_t\cdot v+\nu\nabla u\cdot\nabla v+[(u\cdot\nabla)u]\cdot v+\lambda(v\phi)\cdot\nabla\mu\,dx=\int_{\Omega}h\cdot v\,dx,\\
&\int_{\Omega}\phi_t\psi\,dx+\int_{\Omega}(u\cdot\nabla\phi)\psi\,dx+\gamma\int_{\Omega}\nabla\mu\cdot\nabla\psi\,dx=0,\\
&\int_{\Omega}\nabla\phi\cdot\nabla\theta+f(\phi)\theta\,dx+\int_{\Gamma}\phi_t\theta+\alpha\nabla_\Gamma\phi\cdot\nabla_\Gamma\theta+\beta\phi\theta+g(\phi)\theta\,dS=\int_{\Omega
}\mu\theta\,dx
\end{align*}
for all test functions $v\in V$ and $\psi,$ $\theta\in W=\{w\in H^1(\bar{\Omega},d\sigma):mw=0\}.$ 
\end{definition}
 \section{The well-posedness of weak solutions}
 \def\theequation{3.\arabic{equation}}\makeatother
\setcounter{equation}{0}
In this section, for the sake of completeness, we give the proof of the well-posedness of weak solutions for problem \eqref{1.1}-\eqref{1.7}. Now, we state it as follows.
\begin{theorem}\label{3.1}
Assume that $h\in L^2(\Omega)$ and $(H_1)$-$(H_2)$ hold. Then for any $u_0\in H$ and $\Phi_0=(\phi_0,\theta_0)\in H^1(\bar{\Omega},d\sigma),$ there exists a unique
weak solution $(u(t),\phi(t))$ for problem \eqref{1.1}-\eqref{1.7} such that $m\phi(t)=m\phi_0,$ which depends continuously on the initial data $(u_0,\phi_0,\theta_0)$ with respect to the norm in $H\times H^1(\bar{\Omega},d\sigma).$
\end{theorem}
\textbf{Proof.} We first prove the existence of weak solutions for problem \eqref{1.1}-\eqref{1.7} by the Faedo-Galerkin method (see \cite{cl, tr}).

 Let $A_1=-P\Delta$ is the Stokes operator and $P$ is the Leray-Helmotz projector from $L^2(\Omega)$ onto $H.$ It is well-known that for the eigenvalue problem $A_1\omega=\kappa \omega,$ where there exists a sequences of non-decreasing numbers $\{\kappa_n\}_{n=1}^{\infty}$ and a sequences of functions $\{\omega_n\}_{n=1}^{\infty},$ which are orthonormal and complete in $H$ such that for every $k\geq 1,$ we have
\begin{align*}
A_1\omega_k=\kappa_k \omega_k
\end{align*}
 and
\begin{align*}
\lim_{k\rightarrow+\infty}\kappa_k=+\infty.
\end{align*}
We also introduce the operator $\mathcal{N}$ which is the inverse of the Laplacian operator $-\Delta,$ where $-\Delta$ is endowed with Neumann boundary conditions imposing zero average over the domain $\Omega.$ It is well-known that there exists a sequences of non-decreasing numbers $\{\lambda_n\}_{n=1}^{\infty}$ and a sequences of functions $\{\psi_n\}_{n=1}^{\infty},$ which are orthonormal and complete in $L^2(\Omega)$ such that $\lambda_1=0$ and $\psi_1=1$ as well as for every $k\geq 2,$ we have
\begin{align*}
\mathcal{N}\psi_k=\frac{1}{\lambda_k}\psi_k
\end{align*}
 and
\begin{align*}
 \lim_{k\rightarrow+\infty}\lambda_k=+\infty.
\end{align*}
 For any $n\geq1,$ we introduce two finite-dimensional spaces $W_n=span\{\psi_1, ...,\psi_n\}$ and $H_n=span\{\omega_1, ...,\omega_n\}.$ Let $P_n$ be
the orthogonal projector from $L^2(\Omega)$ to $W_n$ and let $\Pi_n$ be
the orthogonal projector from $H$ to $H_n.$

Consider the approximate solution $(u_n(t),\phi_n(t),\mu_n(t))$ in the form
\begin{align*}
&u_n(t)=\sum_{i=1}^n \beta_i(t)\omega_i,\\
&\phi_n(t)=\sum_{i=1}^n \alpha_i(t)\psi_i,\\
&\mu_n(t)=\sum_{i=1}^n \mu_i(t)\psi_i,
\end{align*}
we obtain $(u_n(t),\phi_n(t))$ from solving the following problem
\begin{equation}\label{3.2}
\begin{cases}
&\int_{\Omega}\frac{\partial u_n}{\partial t}\cdot v+\nu\nabla u_n\cdot\nabla v+[(u_n\cdot\nabla)u_n]\cdot v+\lambda(v\phi_n)\cdot\nabla\mu_n\,dx=\int_{\Omega}h\cdot v\,dx,\\
&\int_{\Omega}\frac{\partial\phi_n}{\partial t}\psi+(u_n\cdot\nabla\phi_n)\psi+\gamma\nabla\mu_n\cdot\nabla\psi\,dx=0,\\
&\int_{\Omega}\nabla\phi_n\cdot\nabla\theta+f(\phi_n)\theta\,dx+\int_{\Gamma}\frac{\partial\phi_n}{\partial t}\theta+\alpha\nabla_\Gamma\phi_n\cdot\nabla_\Gamma\theta+\beta\phi_n\theta+g(\phi_n)\theta\,dS=\int_{\Omega
}\mu_n\theta\,dx,\\
&\int_\Omega u_n(0)\cdot\omega_k\,dx=\int_\Omega u_0\cdot\omega_k\,dx,\, k=1,\cdots,n,\\
&\int_\Omega \phi_n(0)\psi_k\,dx=\int_\Omega \phi_0\psi_k\,dx,\, k=1,\cdots,n,\\
&\int_\Gamma \phi_n(0)\psi_k\,dS=\int_\Gamma \theta_0\psi_k\,dS,\, k=1,\cdots,n,
\end{cases}
\end{equation}
for any $v\in H_n$ and $\psi,$ $\theta\in W_n.$

Repeating the similar argument as in \cite{cl}, we can obtain the local (in time) existence of $(u_n(t),\phi_n(t),\mu_n(t)).$
Next, we will establish some a priori estimates for $(u_n(t),\phi_n(t),\mu_n(t)).$ Let $v=u_n,$ $\psi=\mu_n$ and $\theta=\phi_n$ in equation \eqref{3.2}, we find
 \begin{align}\label{3.3}
\nonumber&\frac{d}{dt}\left(\frac{1}{2}\|u_n\|_{L^2(\Omega)}^2+\frac{\lambda}{2}\|\phi_n\|_{H^1(\bar{\Omega},d\sigma)}^2+\lambda\int_\Omega F(\phi_n)\,dx+\lambda\int_\Gamma G(\phi_n)\,dS\right)+\lambda\|\frac{\partial\phi_n(t)}{\partial t}\|_{L^2
(\Gamma)}^2\\
\nonumber&+\lambda\gamma\|\nabla\mu_n\|_{L^2(\Omega)}^2+\nu\|\nabla u_n\|_{L^2(\Omega)}^2+\|\phi_n\|_{H^1(\bar{\Omega},d\sigma)}^2+\int_\Omega f(\phi_n)\phi_n\,dx+\int_\Gamma g(\phi_n)\phi_n\,dS\\
\nonumber=&\int_\Omega \mu_n\phi_n\,dx-\int_{\Gamma}\frac{\partial\phi_n}{\partial t}\phi_n\,dS+\int_\Omega h(x)u_n\,dx\\
\nonumber\leq&\|\mu_n-m\mu_n\|_{L^2(\Omega)}\|\phi_n\|_{L^2(\Omega)}+m\mu_nm\phi_0|\Omega|+\|\frac{\partial\phi_n}{\partial t}\|_{L^2(\Gamma)}\|\phi_n\|_{L^2(\Gamma)}+\|h\|_{L^2(\Omega)}\|u_n\|_{L^2(\Omega)}\\
\leq&C\|\nabla\mu_n\|_{L^2(\Omega)}\|\phi_n\|_{L^2(\Omega)}+|m\mu_n||m\phi_0||\Omega|+\|\frac{\partial\phi_n}{\partial t}\|_{L^2(\Gamma)}\|\phi_n\|_{L^2(\Gamma)}+\|h\|_{L^2(\Omega)}\|u_n\|_{L^2(\Omega)},
\end{align}
where $F(s)=\int_0^sf(r)\,dr$ and $G(s)=\int_0^sg(r)\,dr,$ respectively, are the primitive function of $f$ and $g.$

Let $\theta=1$ in the third equation of \eqref{3.2}, we obtain
\begin{align}\label{3.4}
\left|\int_\Omega\mu_n\,dx\right|\leq\|\frac{\partial\phi_n(t)}{\partial t}\|_{L^1
(\Gamma)}+\beta\|\phi_n(t)\|_{L^1(\Gamma)}+\|g(\phi_n)\|_{L^1
(\Gamma)}+\|f(\phi_n)\|_{L^1(\Omega)}.
\end{align}
Combining \eqref{1.4}-\eqref{1.7} and \eqref{3.3}-\eqref{3.4} with H\"{o}lder inequality and Young inequality, we find that there exist two positive constants $\delta$ and $\varrho$ such that
 \begin{align}\label{3.5}
\frac{d}{dt}J(u_n,\phi_n)+\lambda\|\frac{\partial\phi_n(t)}{\partial t}\|_{L^2
(\Gamma)}^2+\lambda\gamma\|\nabla\mu_n\|_{L^2(\Omega)}^2+\nu\|\nabla u_n\|_{L^2(\Omega)}^2+\delta J(u_n,\phi_n)\leq\varrho,
\end{align}
where
 \begin{align*}
J(u,\phi)=\|u\|_{L^2(\Omega)}^2+\lambda\|\phi\|_{H^1(\bar{\Omega},d\sigma)}^2+2\lambda\int_\Omega F(\phi)\,dx+2\lambda\int_\Gamma G(\phi)\,dS.
\end{align*}

From the classical Gronwall inequality, we infer
\begin{align}\label{3.6}
\nonumber J(u_n(t),\phi_n(t))\leq &e^{-\delta t}J(u_n(0),\phi_n(0))+\frac{\varrho}{\delta}\\
\leq &e^{-\delta t}J(u_0,\Phi_0)+\frac{\varrho}{\delta}.
\end{align}
From \eqref{1.4}-\eqref{1.7}, we deduce that there exist four positive constants $\delta_1,$ $\delta_2,$ $k_1$ and $k_2$ such that
\begin{align}\label{3.7}
\delta_1\left(\|u\|_{L^2(\Omega)}^2+\lambda\|\phi\|_{H^1(\bar{\Omega},d\sigma)}^2+\|\phi\|_{L^p(\Omega)}^p+\|\phi\|_{L^q(\Gamma)}^q\right)-k_1\leq J(u,\phi)
\leq\delta_2\left(\|u\|_{L^2(\Omega)}^2+\|\phi\|_{H^1(\bar{\Omega},d\sigma)}^2+\|\phi\|_{L^p(\Omega)}^p+\|\phi\|_{L^q(\Gamma)}^q\right)+k_2.
\end{align}
By virtue of \eqref{3.6}-\eqref{3.7}, we obtain
\begin{align}\label{3.8}
\|u_n(t)\|_{L^2(\Omega)}^2+\|\phi_n(t)\|_{H^1(\bar{\Omega},d\sigma)}^2+\|\phi_n(t)\|_{L^p(\Omega)}^p+\|\phi_n(t)\|_{L^q(\Gamma)}^q
\leq\frac{1}{\delta_1}e^{-\delta t}J(u_0,\Phi_0)+\frac{\varrho}{\delta_1\delta}+\frac{k_1}{\delta_1}.
\end{align}
Integrating \eqref{3.5} from $0$ to $t$, we obtain
 \begin{align}\label{3.9}
\lambda\int_0^t\|\frac{\partial\phi_n(s)}{\partial t}\|_{L^2
(\Gamma)}^2\,ds+\lambda\gamma\int_0^t\|\nabla\mu_n(s)\|_{L^2(\Omega)}^2\,ds+\nu\int_0^t\|\nabla u_n(s)\|_{L^2(\Omega)}^2\,ds
\leq\varrho T+J(u_0,\Phi_0)+k_1(1+\delta T),
\end{align}
for any $t\in (0,T]$.

 Due to \eqref{3.4} and \eqref{3.8}-\eqref{3.9}, we find
\begin{align*}
&\{u_n\}_{n=1}^{\infty}\;\textit{is uniformly bounded in}\; L^{\infty}(0,T;H)\cap L^2(0,T;V),\\
&\{\phi_n\}_{n=1}^{\infty}\;\textit{is uniformly bounded in}\;L^{\infty}(0,T;H^1(\bar{\Omega},d\sigma))\cap L^{\infty}(0,T;L^p(\Omega))\cap L^{\infty}(0,T;L^q(\Gamma)),\\
&\{\frac{\partial\phi_n(t)}{\partial t}\}_{n=1}^{\infty}\;\textit{is uniformly bounded in}\;L^2(0,T;L^2(\Gamma)),\\
&\{\mu_n\}_{n=1}^{\infty}\;\textit{is uniformly bounded in}\; L^2(0,T;H^1(\Omega)).
\end{align*}
 Therefore, there exist
 \begin{align*}
&u\in L^{\infty}(0,T;H)\cap L^2(0,T;V),\\
&\phi\in L^{\infty}(0,T;H^1(\bar{\Omega},d\sigma))\cap L^{\infty}(0,T;L^p(\Omega))\cap L^{\infty}(0,T;L^q(\Gamma)),\\
&\frac{\partial\phi}{\partial t}\in L^2(0,T;L^2(\Gamma)),\\
&\chi\in  L^2(0,T;H^1(\Omega))
\end{align*}
  such that we can extract subsequences $\{u_{n_j}\}_{j=1}^{\infty},$ $\{\phi_{n_j}\}_{j=1}^{\infty},$ $\{\frac{\partial\phi_{n_j}}{\partial t}\}_{j=1}^{\infty},$ $\{\mu_{n_j}\}_{j=1}^{\infty}$ of $\{u_n\}_{n=1}^{\infty},$ $\{\phi_n\}_{n=1}^{\infty},$ $\{\frac{\partial\phi_n}{\partial t}\}_{n=1}^{\infty},$ $\{\mu_n\}_{n=1}^{\infty},$ respectively, satisfy
\begin{align*}
&u_{n_j}\rightharpoonup u\;\textit{weakly star in}\; L^{\infty}(0,T;H),\\
&u_{n_j}\rightharpoonup u\;\textit{weakly in}\; L^2(0,T;V),\\
&\phi_{n_j}\rightharpoonup \phi\;\textit{weakly star in}\; L^{\infty}(0,T;H^1(\bar{\Omega},d\sigma)),\\
&\phi_{n_j}\rightharpoonup \phi\;\textit{weakly star in}\;  L^{\infty}(0,T;L^p(\Omega)),\\
&\phi_{n_j}\rightharpoonup \phi\;\textit{weakly star in}\; L^{\infty}(0,T;L^q(\Gamma)),\\
&\frac{\partial\phi_{n_j}(t)}{\partial t}\rightharpoonup \frac{\partial\phi(t)}{\partial t}\;\textit{weakly in}\;L^2(0,T;L^2(\Gamma)),\\
&\mu_{n_j}\rightharpoonup \chi\;\textit{weakly in}\; L^2(0,T;H^1(\Omega)).
\end{align*}
From \eqref{1.4}-\eqref{1.7} and \eqref{3.8}, we obtain
\begin{align}
\label{3.10}&\{f(\phi_n)\}_{n=1}^{\infty}\;\textit{is uniformly bounded in}\;L^\infty(0,T;L^{\frac{p}{p-1}}(\Omega)),\\
\label{3.11}&\{g(\phi_n)\}_{n=1}^{\infty}\;\textit{is uniformly bounded in}\;L^\infty(0,T;L^{\frac{q}{q-1}}(\Gamma)).
\end{align}
We infer from \eqref{1.4}-\eqref{1.7}, \eqref{3.8}-\eqref{3.9}, Sobolev embedding Theorem and Lemma \ref{2.4} that
\begin{align}\label{3.12}
\{\phi_n\}_{n=1}^{\infty}\;\textit{is uniformly bounded in}\;L^2(0,T;H^2(\bar{\Omega},d\sigma)),
\end{align}
entails one can extract a subsequence $\{\phi_{n_j}\}_{j=1}^{\infty}$ of $\{\phi_n\}_{n=1}^{\infty}$ such that
\begin{align*}
&\phi_{n_j}\rightharpoonup \phi\;\textit{weakly in}\; L^2(0,T;H^2(\bar{\Omega},d\sigma)).
\end{align*}
 For any $v\in V,$ set $v_n=\Pi_nv,$ we have
\begin{align*}
|\int_\Omega\frac{\partial u_n}{\partial t}\cdot v\,dx|\leq&\int_\Omega|u_n|^2|\nabla v_n|+\nu|\nabla u_n||\nabla v_n|+|h||v_n|+\lambda|v_n||\phi_n||\nabla\mu_n|\,dx\\
\leq&\|u_n\|_{L^4(\Omega)}^2\|\nabla v_n\|_{L^2(\Omega)}+\nu\|\nabla u_n\|_{L^2(\Omega)}\|\nabla v_n\|_{L^2(\Omega)}\\
&+\|h\|_{L^2(\Omega)}\|v_n\|_{L^2(\Omega)}+\lambda\|v_n\|_{L^{\frac{2p}{p-2}}(\Omega)}\|\phi_n\|_{L^p(\Omega)}\|\nabla\mu_n\|_{L^2(\Omega)}\\
\leq&C\|u_n\|_{L^2(\Omega)}\|\nabla u_n\|_{L^2(\Omega)}\|\nabla v\|_{L^2(\Omega)}+\nu\|\nabla u_n\|_{L^2(\Omega)}\|\nabla v\|_{L^2(\Omega)}\\
&+\frac{1}{\sqrt{\kappa_1}}\|h\|_{L^2(\Omega)}\|\nabla v\|_{L^2(\Omega)}+C\|\nabla v\|_{L^2(\Omega)}\|\phi_n\|_{L^p(\Omega)}\|\nabla\mu_n\|_{L^2(\Omega)},
\end{align*}
entails that
\begin{align*}
\{\frac{\partial u_n}{\partial t}\}_{n=1}^{\infty}\,\,\, \textit{is uniformly bounded in}\,\,\,L^2(0,T;V^{*}).
\end{align*}
 For any $\psi\in H^1(\bar{\Omega},d\sigma),$ set $\psi_n=P_n\psi,$ we have
\begin{align*}
|\int_\Omega\frac{\partial\phi_n}{\partial t}\psi\,dx|\leq&\int_\Omega|u_n\phi_n||\nabla\psi_n|\,dx+\gamma\int_\Omega|\nabla\mu_n||\nabla\psi_n|\,dx\\
\leq&\|u_n\|_{L^{\frac{2p}{p-2}}(\Omega)}\|\phi_n\|_{L^p(\Omega)}\|\nabla\psi_n\|_{L^2(\Omega)}+\gamma\|\nabla\mu_n\|_{L^2(\Omega)}\|\nabla\psi_n\|_{L^2(\Omega)}\\
\leq&C\|\nabla u_n\|_{L^2(\Omega)}\|\phi_n\|_{L^p(\Omega)}\|\nabla\psi\|_{L^2(\Omega)}+\gamma\|\nabla\mu_n\|_{L^2(\Omega)}\|\nabla\psi\|_{L^2(\Omega)},
\end{align*}
which implies that
\begin{align*}
\{\frac{\partial\phi_n}{\partial t}\}_{n=1}^{\infty}\,\,\, \textit{is uniformly bounded in}\,\,\,L^2(0,T;(H^1(\bar{\Omega},d\sigma))^{*}).
\end{align*}
 Therefore, we can extract subsequences $\{\frac{\partial u_{n_j}}{\partial t}\}_{j=1}^{\infty},$ $\{\frac{\partial\phi_{n_j}}{\partial t}\}_{j=1}^{\infty}$ of $\{\frac{\partial u_n}{\partial t}\}_{n=1}^{\infty},$ $\{\frac{\partial\phi_n}{\partial t}\}_{n=1}^{\infty},$ respectively, such that
\begin{align*}
&\frac{\partial u_{n_j}}{\partial t}\rightharpoonup \frac{\partial u}{\partial t}\,\textit{weakly in}\;L^2(0,T;V^{*})\\
&\frac{\partial\phi_{n_j}}{\partial t}\rightharpoonup \frac{\partial\phi}{\partial t}\,\textit{weakly in}\;L^2(0,T;(H^1(\bar{\Omega},d\sigma))^{*}).
\end{align*}
By virtue of the Aubin-Lions compactness theorem, we can extract a further
subsequence (still denote by $\{u_{n_j}\}_{j=1}^{\infty}$ and $\{\phi_{n_j}\}_{j=1}^{\infty}$) such that additionally
\begin{align}
\label{3.13}&u_{n_j}\longrightarrow \phi\;\textit{strongly in}\; L^2(0,T;H),\\
\label{3.14}&\phi_{n_j}\longrightarrow \phi\;\textit{strongly in}\; L^2(0,T;H^1(\bar{\Omega},d\sigma)).
\end{align}
From \eqref{3.10}-\eqref{3.11}, \eqref{3.13}-\eqref{3.14} and Lemma \ref{2.3}, we obtain
\begin{align}
\label{3.15}f(\phi_{n_j})\rightharpoonup f(\phi)\;\textit{weakly in}\;L^{\frac{p-1}{p}}(0,T;L^{\frac{p-1}{p}}(\Omega))\\
\label{3.16}g(\phi_{n_j})\rightharpoonup g(\phi)\;\textit{weakly in}\;L^{\frac{q-1}{q}}(0,T;L^{\frac{q-1}{q}}(\Gamma)).
\end{align}
Hence, we have
\begin{align*}
\chi=-\Delta\phi+f(\phi)=\mu.
\end{align*}
Thanks to
\begin{align*}
&\int_{\Omega_T} v\cdot\left([u_n\cdot\nabla]u_n-[u\cdot\nabla]u\right)\,dx=\int_{\Omega_T} v\cdot([(u_n-u)\cdot\nabla]u_n)\,dx+\int_{\Omega_T} v\cdot([u\cdot\nabla](u_n-u))\,dx,\\
&\int_{\Omega_T} v\cdot(\phi_n\nabla\mu_n-\phi\nabla\mu)\,dx=\int_{\Omega_T} v\cdot\nabla\mu_n(\phi_n-\phi)\,dx+\int_{\Omega_T} v\phi\cdot(\nabla\mu_n-\nabla\mu)\,dx
\end{align*}
for any $v\in V$ and
\begin{align*}
\int_{\Omega_T}(u_n\cdot\nabla\phi_n-u\cdot\nabla\phi)\psi\,dx=-\int_{\Omega_T}\phi(u_n-u)\cdot\nabla\psi\,dx-\int_{\Omega_T} (\phi_n-\phi)u_n\cdot\nabla\psi\,dx
\end{align*}
for any $\psi\in H^1(\bar{\Omega},d\sigma),$ which imply
\begin{align*}
&(u_n\cdot\nabla)u_n\rightharpoonup (u\cdot\nabla)u\,\textit{weakly in}\;L^2(0,T;V^{*}),\\
&\phi_n\nabla\mu_n\rightharpoonup\phi\nabla\mu\,\textit{weakly in}\;L^2(0,T;V^*),\\
&u_n\cdot\nabla\phi_n\rightharpoonup u\cdot\nabla\phi\,\textit{weakly in}\;L^2(0,T;(H^1(\bar{\Omega},d\sigma))^{*}),\\
\end{align*}
Therefore, a weak solution $(u,\phi)$ for problem \eqref{1.1}-\eqref{1.7} has been proved, and we infer $u(t)\in\mathcal{C}(\mathbb{R}^+;H)$ and $\phi(t)\in\mathcal{C}(\mathbb{R}^+;H^1(\bar{\Omega},d\sigma))$ from lemma \ref{2.5}.

Finally, we prove the uniqueness and the continuous dependence on the initial data of the solutions. Let $(u_1,\phi_1,p_1)$, $(u_2,\phi_2,p_2)$
be two solutions for problem \eqref{1.1}-\eqref{1.7} with the initial data
$(u_{1_0},\phi_{1_0},\theta_{1_0}),$ $(u_{2_0},\phi_{2_0},\theta_{2_0}),$ respectively, and $m\phi_{1_0}=m\phi_{2_0}.$ Let $u=u_1-u_2,$ $\phi=\phi_1-\phi_2,$ $p=p_1-p_2,$ then $(u,\phi,p)$ satisfies the following equations
\begin{equation}\label{3.17}
\begin{cases}
&\frac{\partial u}{\partial t}-\nu\Delta u+u\cdot\nabla u_2+u_1\cdot\nabla u+\nabla p =-\lambda\phi_1\nabla\mu-\lambda\phi\nabla\mu_2,\,\,\,(x,t)\in\Omega\times\mathbb{R}^+,\\
&\nabla\cdot u=0,\,\,\,(x,t)\in\Omega\times\mathbb{R}^+,\\
&\frac{\partial \phi}{\partial t}+u\cdot\nabla\phi_1+u_2\cdot\nabla\phi-\gamma\Delta\mu=0,\,\,\,(x,t)\in\Omega\times\mathbb{R}^+,\\
&\mu=\mu_1-\mu_2=-\Delta\phi+f(\phi_1)-f(\phi_2),\,\,\,(x,t)\in\Omega\times\mathbb{R}^+.
\end{cases}
\end{equation}
Equation \eqref{3.17} is subject to the following boundary conditions
\begin{equation}\label{3.18}
\begin{cases}
&u(x,t)=0,\,\,\,(x,t)\in\Gamma\times\mathbb{R}^+,\\
&\frac{\partial \mu}{\partial \vec{n}}=0,\,\,\,(x,t)\in\Gamma\times\mathbb{R}^+,\\
&\frac{\partial \phi}{\partial t}-\alpha\Delta_\Gamma\phi+\frac{\partial \phi}{\partial \vec{n}}+\beta\phi+g(\phi_1)-g(\phi_2)=0,\,\,\,(x,t)\in\Gamma\times\mathbb{R}^+
\end{cases}
\end{equation}
 and initial conditions
\begin{equation}\label{3.19}
\begin{cases}
 &u(x,0)=u_{1_0}-u_{2_0},\,\,\,x\in\Omega,\\
 &\phi(x,0)=\phi_{1_0}-\phi_{2_0},\,\,\,x\in\Omega,\\
 &\phi(x,0)=\theta_{1_0}-\theta_{1_0},\,\,\,x\in\Gamma.
\end{cases}
\end{equation}
Multiplying the first equation and the third equation of \eqref{3.17} by $u,$ $-\lambda\Delta\phi,$ respectively, and integrating by parts, we find
\begin{align}\label{3.20}
\nonumber&\frac{1}{2}\frac{d}{dt}(\|u(t)\|_{L^2(\Omega)}^2+\lambda\|\phi(t)\|_{H^1(\bar{\Omega},d\sigma)}^2)+\lambda\|\phi_t(t)\|_{L^2(\Gamma)}^2+\nu\|\nabla u\|_{L^2(\Omega)}^2+\lambda\gamma\|\nabla\Delta\phi\|_{L^2(\Omega)}^2\\
\nonumber=&\lambda\gamma\int_{\Omega}\nabla(f(\phi_1)-f(\phi_2))\cdot\nabla\Delta\phi\,dx-\lambda\int_{\Gamma}(g(\phi_1)-g(\phi_2))\phi_t\,dS-\lambda\int_{\Omega}(u_2\phi)\cdot\nabla\Delta\phi\,dx\\
\nonumber&+\lambda\int_{\Omega}(u\cdot\nabla\phi_1)(f(\phi_1)-f(\phi_2))+(u\phi)\cdot\nabla\mu_2\,dx-\int_{\Omega}[(u\cdot\nabla)u_2]\cdot u\,dx\\
\nonumber\leq&\lambda\gamma\int_{\Omega}\nabla(f(\phi_1)-f(\phi_2))\cdot\nabla\Delta\phi\,dx+\lambda\|g(\phi_1)-g(\phi_2)\|_{L^2(\Gamma)}\|\phi_t\|_{L^2(\Gamma)}\\
\nonumber&+\lambda\|u\|_{L^4(\Omega)}\|\nabla\phi_1\|_{L^2(\Omega)}\|f(\phi_1)-f(\phi_2)\|_{L^4(\Omega)}+\lambda\|u_2\|_{L^4(\Omega)}\|\phi\|_{L^4(\Omega)}\|\nabla\Delta\phi\|_{L^2(\Omega)}\\
&+\lambda\|u\|_{L^4(\Omega)}\|\phi\|_{L^4(\Omega)}\|\nabla\mu_2\|_{L^2(\Omega)}+\|u\|_{L^4(\Omega)}^2\|\nabla u_2\|_{L^2(\Omega)}.
\end{align}
Due to
\begin{align}\label{3.21}
\nonumber\left|\int_{\Omega}\nabla(f(\phi_1)-f(\phi_2))\cdot\nabla\Delta\phi\,dx\right|
\leq&\left|\int_{\Omega}(f'(\phi_1)-f'(\phi_2))\nabla\phi_1\cdot\nabla\Delta\phi\,dx\right|+\left|\int_{\Omega}f'(\phi_2)\nabla\phi\cdot\nabla\Delta\phi\,dx\right|\\
\nonumber\leq&C\|\nabla\phi_1\|_{L^6(\Omega)}(1+\|\phi_1\|_{L^{6(p-3)}(\Omega)}^{p-3}+\|\phi_2\|_{L^{6(p-3)}(\Omega)}^{p-3})\|\phi\|_{L^6(\Omega)}\|\nabla\Delta\phi\|_2\\
\nonumber&+C(1+\|\phi_2\|_{L^{4(p-2)}(\Omega)}^{p-2})\|\nabla\phi\|_{L^4(\Omega)}\|\nabla\Delta\phi\|_2\\
\nonumber\leq&C\|\phi_1\|_{H^2(\bar{\Omega},d\sigma)}(1+\|\phi_1\|_{H^1(\bar{\Omega},d\sigma)}^{p-3}+\|\phi_2\|_{H^1(\bar{\Omega},d\sigma)}^{p-3})\|\phi\|_{H^1(\bar{\Omega},d\sigma)}\|\nabla\Delta\phi\|_2\\
\nonumber&+C(1+\|\phi_2\|_{H^1(\bar{\Omega},d\sigma)}^{p-2})\|\nabla\phi\|_{L^4(\Omega)}\|\nabla\Delta\phi\|_2\\
\nonumber\leq&C\|\phi_1\|_{H^2(\bar{\Omega},d\sigma)}(1+\|\phi_1\|_{H^1(\bar{\Omega},d\sigma)}^{p-3}+\|\phi_2\|_{H^1(\bar{\Omega},d\sigma)}^{p-3})\|\phi\|_{H^1(\bar{\Omega},d\sigma)}\|\nabla\Delta\phi\|_2\\
\nonumber&+C(1+\|\phi_2\|_{H^1(\bar{\Omega},d\sigma)}^{p-2})(\|\phi\|_{L^2(\Omega)}+\|\phi\|_{L^2(\Omega)}^{\frac{1}{2}}\|\nabla\Delta\phi\|_2^{\frac{1}{2}})\|\nabla\Delta\phi\|_2\\
\nonumber\leq&C\|\phi_1\|_{H^2(\bar{\Omega},d\sigma)}(1+\|\phi_1\|_{H^1(\bar{\Omega},d\sigma)}^{p-3}+\|\phi_2\|_{H^1(\bar{\Omega},d\sigma)}^{p-3})\|\phi\|_{H^1(\bar{\Omega},d\sigma)}\|\nabla\Delta\phi\|_2\\
&+C(1+\|\phi_2\|_{H^1(\bar{\Omega},d\sigma)}^{p-2})(\|\phi\|_{H^1(\bar{\Omega},d\sigma)}+\|\phi\|_{H^1(\bar{\Omega},d\sigma)}^{\frac{1}{2}}\|\nabla\Delta\phi\|_2^{\frac{1}{2}})\|\nabla\Delta\phi\|_2,
\end{align}
where we use the following Gagliardo-Nirenberg inequality:
\begin{align*}
\|\nabla\phi\|_{L^4(\Omega)}\leq C\|\nabla\Delta\phi\|_{L^2(\Omega)}^{\frac{1}{2}}\|\phi\|_{L^2(\Omega)}^{\frac{1}{2}}+C_2\|\phi\|_{L^2(\Omega)},
\end{align*}
\begin{align}\label{3.22}
\nonumber\|f(\phi_1)-f(\phi_2)\|_{L^4(\Omega)}\leq&C(1+\|\phi_1\|_{L^{8(p-2)}(\Omega)}^{p-2}+\|\phi_2\|_{L^{8(p-2)}(\Omega)}^{p-2})\|\phi\|_{L^8(\Omega)}\\
\leq&C(1+\|\phi_1\|_{H^1(\bar{\Omega},d\sigma)}^{p-2}+\|\phi_2\|_{H^1(\bar{\Omega},d\sigma)}^{p-2})\|\phi\|_{H^1(\bar{\Omega},d\sigma)}
\end{align}
and
\begin{align}\label{3.23}
\nonumber\|g(\phi_1)-g(\phi_2)\|_{L^2(\Gamma)}\leq&C(1+\|\phi_1\|_{L^{4(q-2)}(\Gamma)}^{q-2}+\|\phi_2\|_{L^{4(q-2)}(\Gamma)}^{q-2})\|\phi\|_{L^4(\Gamma)}\\
\leq&C(1+\|\phi_1\|_{H^1(\bar{\Omega},d\sigma)}^{q-2}+\|\phi_2\|_{H^1(\bar{\Omega},d\sigma)}^{q-2})\|\phi\|_{H^1(\bar{\Omega},d\sigma)},
\end{align}
we infer from \eqref{3.20}-\eqref{3.23} and Young inequality that
\begin{align}\label{3.24}
\frac{d}{dt}(\|u(t)\|_{L^2(\Omega)}^2+\lambda\|\phi(t)\|_{H^1(\bar{\Omega},d\sigma)}^2)
\leq\mathbb{L}(t)(\|u(t)\|_{L^2(\Omega)}^2+\lambda\|\phi(t)\|_{H^1(\bar{\Omega},d\sigma)}^2),
\end{align}
where
\begin{align}\label{3.35}
\mathbb{L}(t)=C(1+\|\phi_1\|_{H^2(\bar{\Omega},d\sigma)}^2+\|\phi_2\|_{H^2(\bar{\Omega},d\sigma)}^2+\|\nabla u_2\|_{L^2(\Omega)}^2+\|\nabla\mu_2\|_{L^2(\Omega)}^2).
\end{align}
Therefore, we conclude from \eqref{3.9}, \eqref{3.12} and \eqref{3.35} that
\begin{align*}
\int_0^T\mathbb{L}(s)\,ds=\mathcal{M}(T)<\infty.
\end{align*}
From the classical Gronwall inequality, we obtain
\begin{align*}
\|\phi(t)\|_{H^1(\bar{\Omega},d\sigma)}^2+\|u(t)\|_{L^2(\Omega)}^2
\leq\max\{1,\alpha,\beta\}\left(\|u_{1_0}-u_{2_0}\|_{L^2(\Omega)}^2+\|\nabla\phi_{1_0}-\nabla\phi_{2_0}\|_{L^2(\Omega)}^2+\|\theta_{1_0}-\theta_{2_0}\|_{H^1(\Gamma)}^2\right)e^{\mathcal{M}(T)}.
\end{align*}
Therefore, $(u_1(x,t),\phi_1(x,t))=(u_2(x,t),\phi_2(x,t))$ a.e. in $\overline{\Omega_T}$,
if $u_{1_0}(x)=u_{2_0}(x),$ $\phi_{1_0}(x)=\phi_{2_0}(x)$ in $\Omega$ and $\theta_{1_0}(x)=\theta_{2_0}(x)$ in $\Gamma,$ and $(u(x,t),\phi(x,t))$ depends
continuously on the initial data $(u_0,\phi_0,\theta_0)$ with respect to the norm in $H\times H^1(\bar{\Omega},d\sigma).$
The proof of Theorem \ref{3.1} is completed.\\
\hfill\qed

\begin{corollary}\label{3.36}
Assume that $h\in L^2(\Omega),$ $(u_{0m},\phi_{0m},\theta_{0m})\rightharpoonup(u_0,\phi_0,\theta_0)$ in $H\times H^1(\bar{\Omega},d\sigma)$ and $(H_1)$-$(H_2)$ hold, let $(u_m(t),\phi_m(t))$ be a sequence of weak solution for problem \eqref{1.1}-\eqref{1.7} such that $(u_m(0),\phi_m(0))=(u_{0m},\phi_{0m},\theta_{0m}).$ For any $T>0,$ if there exists a subsequence converging ($\ast$-) weakly in spaces $\{(u,\phi)\in L^{\infty}(0,T;H\times H^1(\bar{\Omega},d\sigma))\cap L^2(0,T; V\times H^2(\bar{\Omega},d\sigma)):(u_t,\phi_t)\in L^1(0,T;(V\times H^1(\bar{\Omega},d\sigma))^*)\}$ to a certain function $(u(t),\phi(t)).$
Then $(u(t),\phi(t))$ is a weak solution on $[0,T]$ with $(u(0),\phi(0))=(u_0,\phi_0,\theta_0).$
\end{corollary}

  For every fixed $I\in\mathbb{R},$ let $V_I=\{\phi\in H^1(\bar{\Omega},d\sigma):m\phi=I\},$ by Theorem \ref{3.1}, we can define the operator semigroup $\{S_I(t)\}_{t\geq 0}$ in $H\times V_I$ by
\begin{align*}
S_I(t)(u_0,\phi_0,\theta_0)=(u(t),\phi(t))=(u(t;(u_0,\phi_0,\theta_0)),\phi(t;(u_0,\phi_0,\theta_0)))
\end{align*}
for all $t\geq0,$ which is $(H\times V_I,H\times V_I)$-continuous, where $(u(t),\phi(t))$ is the solution of problem \eqref{1.1}-\eqref{1.7} with $(u(x,0),\phi(x,0))=(u_0,\phi_0,\theta_0)\in H\times V_I.$

\section{The existence of global attractors}
\def\theequation{4.\arabic{equation}}\makeatother
\setcounter{equation}{0}
\subsection{The existence of a global attractor in $X_\ell$}
In this subsection, we will consider the existence of global attractors for problem \eqref{1.1}-\eqref{1.7} by using the $\ell$-trajectory method. From Theorem \ref{3.1}, we know that the solution $(u(t),\phi(t))$ of problem \eqref{1.1}-\eqref{1.7} with initial data $(u_0,\phi_0,\theta_0)$ in $H\times V_I$ is unique. Therefore, for any $\ell>0$ and any $(u_0,\phi_0,\theta_0)\in H\times V_I,$ there is only one solution defined on the time interval $[0,\ell]$ starting from the initial data $(u_0,\phi_0,\theta_0)\in H\times V_I,$ for the sake of simplicity, which is denoted by $\chi(\tau,(u_0,\phi_0,\theta_0)).$ Denote by $X_\ell$ the set of all the solution trajectories defined on the time interval $[0,\ell]$ equipped with the topology of $L^2(0,\ell;H\times V_I).$ Since $X_\ell\subset\mathcal{C}([0,\ell];H\times V_I),$ it makes sense to talk about the point values of
trajectories. On the other hand, it is not clear whether $X_\ell$ is closed in $L^2(0,\ell;H\times V_I)$
and hence $X_\ell$ in general is not a complete metric space. In what follows, we first give the definition of some operators.

For any $t\in [0,1],$  we define the mapping $e_t: X_\ell\rightarrow H\times V_I$ by
\begin{align*}
e_t(\chi)=\chi(t\ell)
\end{align*}
for any $\chi\in X_\ell.$

The mapping $b: H\times V_I\rightarrow X_\ell$ is given by
\begin{align*}
b((u_0,\phi_0,\theta_0))=S_I(\tau,(u_0,\phi_0,\theta_0)),\,\,\,\tau\in[0,\ell]
\end{align*}
for any $(u_0,\phi_0,\theta_0)\in H\times V_I$ and we define the operators $L_t:X_{\ell}\rightarrow X_{\ell}$ by the relation
\begin{align*}
L_tb((u_0,\phi_0,\theta_0))=(u,\phi)(t+\tau,(u_0,\phi_0,\theta_0))=S_I(t+\tau)(u_0,\phi_0,\theta_0),\,\,\,\tau\in[0,\ell]
\end{align*}
for any $(u_0,\phi_0,\theta_0)\in H\times V_I,$ where $(u,\phi)$ is the unique solution of problem \eqref{1.1}-\eqref{1.7} with initial data $(u_0,\phi_0,\theta_0),$ we can easily prove the operators $\{L_t\}_{t\geq0}$ is a semigroup on $X_\ell.$

Next, we will carry out some a priori estimates to obtain the existence of absorbing sets for problem \eqref{1.1}-\eqref{1.7}.
\begin{theorem}\label{4.1.1}
Assume that $h\in L^2(\Omega)$ and $(H_1)$-$(H_2)$ hold. Then there exists a positive constant $\rho_1$ satisfying for any bounded subset $B\subset H\times V_I,$ there exists a time $\tau_1=\tau_1(B)>0$ such that for any weak solutions of problem \eqref{1.1}-\eqref{1.7} with initial data $(u_0,\phi_0,\theta_0)\in B,$ we have
\begin{align*}
\|u(t)\|_{L^2(\Omega)}^2+\lambda\|\phi(t)\|_{H^1(\bar{\Omega},d\sigma)}^2\leq\rho_1
\end{align*}
and
\begin{align*}
\int_0^{\ell}\|u(t+s)\|_{L^2(\Omega)}^2+\lambda\|\phi(t+s)\|_{H^1(\bar{\Omega},d\sigma)}^2\,ds\leq\rho_1
\end{align*}
for any $t\geq \tau_1.$
\end{theorem}
\textbf{Proof.} From \eqref{3.6}, we infer that for any bounded subset $B\subset H\times V_I,$ there exists a time $\tau_0=\tau_0(B)>0$ such that
\begin{align*}
J(u(t),\phi(t))\leq \delta_1+\frac{\varrho}{\delta}
\end{align*}
for any $t\geq \tau_0,$ which implies that
\begin{align*}
\|u(t)\|_{L^2(\Omega)}^2+\lambda\|\phi(t)\|_{H^1(\bar{\Omega},d\sigma)}^2\leq 1+\frac{\varrho}{\delta_1\delta}+\frac{k_1}{\delta_1}
\end{align*}
for any $t\geq \tau_0.$

From \eqref{3.5}, we deduce
 \begin{align}\label{4.1.2}
\frac{d}{dt}J(u,\phi)+\lambda\|\phi_t(t)\|_{L^2
(\Gamma)}^2+\lambda\gamma\|\nabla\mu\|_{L^2(\Omega)}^2+\nu\|\nabla u\|_{L^2(\Omega)}^2+\delta J(u,\phi)\leq\varrho.
\end{align}
Integrating \eqref{4.1.2} from $0$ to $\ell$ and combining \eqref{3.7}, we obtain
\begin{align}\label{4.1.3}
\lambda\int_0^\ell\|\phi_t(r)\|_{L^2
(\Gamma)}^2\,dr+\delta\int_0^\ell J(u(r),\phi(r))\,dr\leq\varrho\ell+J(u(0),\phi(0))+k_1.
\end{align}
Integrating \eqref{4.1.3} from $r$ to $t+r$ and integrating the resulting inequality with respect to $r$ over $(0,\ell),$ we obtain
\begin{align*}
\int_0^\ell J(u(t+r),\phi(t+r))\,dr \leq &e^{-\delta t}\int_0^\ell J(u(r),\phi(r))\,dr
+\ell\frac{\varrho}{\delta}(1-e^{-\delta t})\\
\leq&e^{-\delta t}\frac{1}{\delta}(\varrho\ell+J(u(0),\phi(0))+k_1)+\ell\frac{\varrho}{\delta},
\end{align*}
which implies that
\begin{align*}
\int_0^{\ell}\|u(t+r)\|_{L^2(\Omega)}^2+\lambda\|\phi(t+r)\|_{H^1(\bar{\Omega},d\sigma)}^2\,dr \leq e^{-\delta t}\frac{1}{\delta\delta_1}(\varrho\ell+J(u(0),\phi(0))+k_1)+\ell\frac{\varrho}{\delta\delta_1}+\frac{k_1}{\delta_1}.
\end{align*}
Therefore, for any bounded subset $B\subset H\times V_I,$ there exists a time $\tau_1=\tau_1(B)>\tau_0$ such that
\begin{align}\label{4.1.4}
\int_0^{\ell}\|u(t+r)\|_{L^2(\Omega)}^2+\lambda\|\phi(t+r)\|_{H^1(\bar{\Omega},d\sigma)}^2\,dr \leq 1+\ell\frac{\varrho}{\delta\delta_1}+\frac{k_1}{\delta_1}
\end{align}
for any $t\geq \tau_1.$\\
\qed\hfill

Let
\begin{align*}
B_0=\left\{(u,\phi)\in H\times V_I:\|u\|_{L^2(\Omega)}^2+\lambda\|\phi\|_{H^1(\bar{\Omega},d\sigma)}^2\leq\rho_1\right\},
\end{align*}
we infer from Theorem \ref{4.1.1} that there exists a time $t_0=t_0(B_0)\geq 0$ such that for any $t\geq t_0,$ we have
\begin{align*}
S_I(t)B_0\subset B_0.
\end{align*}
Define
\begin{align*}
B_1=\overline{\bigcup\limits_{t\in [0,t_0]}S_I(t)B_0}^{H\times V_I}
\end{align*}
and
\begin{align*}
B_0^\ell=\{\chi\in X_\ell:e_0(\chi)\in B_1\},
\end{align*}
from the continuity of $S_I(t),$ \eqref{3.6} and Theorem \ref{4.1.1}, we deduce
\begin{align*}
S_I(t)B_1\subset B_1
\end{align*}
and
\begin{align*}
L_tB_0^\ell\subset B_0^\ell
\end{align*}
for any $t\geq 0$ as well as $B_1$ is a bounded subset of $H\times V_I.$

From Theorem \ref{4.1.1}, we immediately obtain the following result.
\begin{corollary}\label{4.1.5}
Assume that $h\in L^2(\Omega)$ and $(H_1)$-$(H_2)$ hold. Then for any bounded subset $B^\ell\subset X_\ell,$ there exists a time $t_1=t_1(B^\ell)>0$ such that for any weak solutions of problem \eqref{1.1}-\eqref{1.7} with short trajectory $\chi\in B^\ell,$ we have
\begin{align*}
\|u(t)\|_{L^2(\Omega)}^2+\lambda\|\phi(t)\|_{H^1(\bar{\Omega},d\sigma)}^2\leq\rho_1
\end{align*}
and
\begin{align*}
\int_0^{\ell}\|u(t+s)\|_{L^2(\Omega)}^2+\lambda\|\phi(t+s)\|_{H^1(\bar{\Omega},d\sigma)}^2\,ds\leq\rho_1
\end{align*}
for any $t\geq t_1.$
\end{corollary}
In what follows, we prove the existence of a compact absorbing set in $X_\ell$ of the semigroup $\{L_t\}_{t\geq 0}.$
\begin{theorem}\label{4.1.6}
Assume that $h\in L^2(\Omega)$ and $(H_1)$-$(H_2)$ hold. Then there exists a positive constant $\rho_2$ satisfying for the subset $B_0^\ell,$ there exists a time $\tau_2=\tau_2(B_0^\ell)>0$ such that for any weak solutions of problem \eqref{1.1}-\eqref{1.7} with short trajectory $\chi\in B_0^\ell,$ we have
\begin{align*}
\int_0^{\ell}\|\nabla u(t+r)\|_{L^2(\Omega)}^2+\lambda\|\phi(t+r)\|_{H^2(\bar{\Omega},d\sigma)}^2\,dr+\left(\int_0^{\ell}\|u_t(t+r)\|_{V^*}+\|\phi_t(t+r)\|_{V_I^*}\,dr\right)^2\leq\rho_2
\end{align*}
for any $t\geq \tau_2.$
\end{theorem}
\textbf{Proof.} From the proof of Theorem \ref{4.1.1} and Corollary \ref{4.1.6}, we know that there exists a $t_0=t_0(B_0^\ell)$ such that
\begin{align}\label{4.1.7}
\|u(t)\|_{L^2(\Omega)}^2+\lambda\|\phi(t)\|_{H^1(\bar{\Omega},d\sigma)}^2+\int_0^\ell J(u(t+r),\phi(t+r))\,dr\leq 2+\ell\frac{\varrho}{\delta}+\frac{\varrho}{\delta_1\delta}+\frac{k_1}{\delta_1}
\end{align}
for any $t\geq t_0.$

Integrating \eqref{4.1.2} between $t-s$ and $t+\ell$ with $t\geq t_0+\frac{\ell}{2},$ $s\in(0,\frac{\ell}{2}),$ we obtain
\begin{align}\label{4.1.8}
\nonumber&\lambda\int_0^\ell\|\phi_t(t+r)\|_{L^2
(\Gamma)}^2\,dr+\lambda\gamma\int_0^\ell\|\nabla\mu(t+r)\|_{L^2(\Omega)}^2\,dr+\nu\int_0^\ell\|\nabla u(t+r)\|_{L^2(\Omega)}^2\\
\leq &J(u(t-s),\phi(t-s))+\varrho(l+s)+k_1+\delta k_1\ell.
\end{align}
After integrating \eqref{4.1.8} with respect to $s$ over $(0,\frac{\ell}{2})$ and combining \eqref{4.1.7}, we have
 \begin{align}\label{4.1.9}
\lambda\int_0^\ell\|\phi_t(t+r)\|_{L^2
(\Gamma)}^2\,dr+\lambda\gamma\int_0^\ell\|\nabla\mu(t+r)\|_{L^2(\Omega)}^2\,dr+\nu\int_0^\ell\|\nabla u(t+r)\|_{L^2(\Omega)}^2\leq\varrho_1
 \end{align}
for any $t\geq t_0+\frac{\ell}{2}.$

It follows from Lemma \ref{2.4} and \eqref{4.1.7}, \eqref{4.1.9} that
\begin{align}\label{4.1.10}
\int_0^{\ell}\|\nabla u(t+s)\|_{L^2(\Omega)}^2+\lambda\|\phi(t+s)\|_{H^2(\bar{\Omega},d\sigma)}^2\,ds\leq\rho_2
\end{align}
for any $t\geq t_0+\frac{\ell}{2}.$

From the proof of Theorem \ref{3.1}, we conclude
\begin{align}\label{4.1.11}
\|u_t\|_{V^*}\leq C\|u\|_{L^2(\Omega)}\|\nabla u\|_{L^2(\Omega)}+\nu\|\nabla u\|_{L^2(\Omega)}+\frac{1}{\sqrt{\kappa_1}}\|h\|_{L^2(\Omega)}+C\|\phi\|_{V_I}\|\nabla\mu\|_{L^2(\Omega)}
\end{align}
and
\begin{align}\label{4.1.12}
\|\phi_t\|_{V_I^*}\leq C\|\nabla u\|_{L^2(\Omega)}\|\phi\|_{V_I}+\gamma\|\nabla\mu\|_{L^2(\Omega)}.
\end{align}
Integrating \eqref{4.1.11}-\eqref{4.1.12} over $(t,t+\ell)$ and combining \eqref{4.1.9}-\eqref{4.1.10} with H\"{o}lder inequality, we obtain
\begin{align}\label{4.1.13}
\int_0^{\ell}\|u_t(t+r)\|_{V^*}+\|\phi_t(t+r)\|_{V_I^*}\,dr\leq\varrho_3
\end{align}
for any $t\geq t_0+\frac{\ell}{2}.$\\
\qed\hfill

Let
\begin{align*}
Y=\left\{\chi\in X_\ell:\chi\in L^2(0,\ell; V\times H^2(\bar{\Omega},d\sigma),\chi_t\in L^1(0,\ell; V^*\times (H^1(\bar{\Omega},d\sigma))^*\right\}
\end{align*}
equipped with the following norm
\begin{align*}
\|\chi\|_Y=\left\{\int_0^\ell\|\chi(r)\|_{V\times H^2(\bar{\Omega},d\sigma)}^2dr+\left(\int_0^\ell\|\chi_t(r)\|_{V^*\times (H^1(\bar{\Omega},d\sigma))^*}\,dr\right)^2\right\}^{\frac{1}{2}}.
\end{align*}
Define
\begin{align*}
B_1^\ell=\left\{\chi\in X_\ell:\|\chi\|_Y^2\leq\rho_2\right\}.
\end{align*}
From Theorem \ref{4.1.1} and Theorem \ref{4.1.6}, we know that $L_t B_0^\ell\subset B_0^\ell$ for any $t\geq 0$ as well as $L_t B_0^\ell\subset B_1^\ell$ for any $t\geq \tau_2.$
\begin{lemma}
Assume that $h\in L^2(\Omega)$ and $(H_1)$-$(H_2)$ hold. Then $\overline{L_t B_0^\ell}^{L^2(0,\ell;H\times V_I)}\subset B_0^\ell$ for any $t\geq 0.$
\end{lemma}
\textbf{Proof.} Thanks to $L_t B_0^\ell\subset B_0^\ell$ for any $t\geq 0,$ it is enough to prove that
\begin{align*}
\overline{B_0^\ell}^{L^2(0,\ell;H\times V_I)}\subset B_0^\ell.
\end{align*}
For any $\chi_0\in \overline{B_0^\ell}^{L^2(0,\ell;H\times V_I)},$ there exists a sequence of trajectories $\chi_n\in B_0^\ell$ such that $\chi_n\rightarrow\chi_0$ in $L^2(0,\ell;H\times V_I),$ which implies that $e_t(\chi_n)\rightarrow e_t(\chi_0)$ in $H\times V_I$ for almost all $t\in [0,1].$ Since $e_0(\chi_n)\in B_1$ for any $n\in\mathbb{N},$ there exists a subsequence $\{e_0(\chi_{n_j})\}_{j=1}^{\infty}$ of $\{e_0(\chi_n)\}_{n=1}^{\infty}$ and $(u_0,\phi_0,\theta_0)\in H\times V_I$ such that $e_0(\chi_{n_j})\rightharpoonup (u_0,\phi_0,\theta_0)$ in $H\times V_I.$ From the proof of the existence of weak solutions for problem \eqref{1.1}-\eqref{1.7}, we deduce that for any $T>0,$ there exists a subsequence converging ($\ast$-) weakly in spaces $\{(u,\phi)\in L^{\infty}(0,T;H\times H^1(\bar{\Omega},d\sigma))\cap L^2(0,T; V\times H^2(\bar{\Omega},d\sigma)):(u_t,\phi_t)\in L^1(0,T;((V\times H^1(\bar{\Omega},d\sigma))*)\}$ to a certain function $(u(t),\phi(t))$ with $(u(0),\phi(0))=(u_0,\phi_0,\theta_0).$ Therefore, we obtain $\chi_0\in X_\ell$ from Corollary \ref{3.36}. It remains to show that $e_0(\chi)\in B_1.$ Since $B_1$ is closed, $e_t(\chi_0)\in B_1$ for almost all $t\in [0,1].$ In particular, $e_{t_n}(\chi_0)\in B_1$ for any sequence $t_n$ with $t_n\rightarrow 0.$ From the continuity of $\chi_0:[0,\ell]\rightarrow H\times V_I$ and the closedness of $B_1,$ we deduce that $e_0(\chi_0)\in B_1.$ Therefore, we obtain $\chi_0\in B_0^\ell.$\\
\qed\hfill

\begin{lemma}\label{4.1.14}
Assume that $h\in L^2(\Omega)$ and $(H_1)$-$(H_2)$ hold. Then the mapping $L_t:X_{\ell}\rightarrow X_{\ell}$ is locally Lipschitz continuous on $B_1^\ell$ for all $t\geq0.$
\end{lemma}
\textbf{Proof.} For any fixed $t>0,$ let $(u_1,\phi_1,p_1)$, $(u_2,\phi_2,p_2)$
be two solutions for problem \eqref{1.1}-\eqref{1.7} with the initial data
$(u_{1_0},\phi_{1_0},\theta_{1_0}),$ $(u_{2_0},\phi_{2_0},\theta_{2_0}),$ respectively, and $m\phi_{1_0}=m\phi_{2_0}.$ Let $u=u_1-u_2,$ $\phi=\phi_1-\phi_2,$ $p=p_1-p_2,$ from the proof of Theorem \ref{3.1}, we conclude
\begin{align}\label{4.1.15}
\frac{d}{dt}(\|u(t)\|_{L^2(\Omega)}^2+\lambda\|\phi(t)\|_{H^1(\bar{\Omega},d\sigma)}^2)
\leq\mathbb{L}(t)(\|u(t)\|_{L^2(\Omega)}^2+\lambda\|\phi(t)\|_{H^1(\bar{\Omega},d\sigma)}^2),
\end{align}
where
\begin{align*}
\mathbb{L}(t)=C(1+\|\phi_1\|_{H^2(\bar{\Omega},d\sigma)}^2+\|\phi_2\|_{H^2(\bar{\Omega},d\sigma)}^2+\|\nabla u_2\|_{L^2(\Omega)}^2+\|\nabla\mu_2\|_{L^2(\Omega)}^2).
\end{align*}
Let $s\in(0,\ell)$ and integrating \eqref{4.1.15} from $s$ to $t+s,$ we
obtain
\begin{align}\label{4.1.16}
\nonumber&\|u(t+s)\|_{L^2(\Omega)}^2+\lambda\|\phi(t+s)\|_{H^1(\bar{\Omega},d\sigma)}^2\\
\leq&\int_s^{t+s}\mathbb{L}(r)(\|u(r)\|_{L^2(\Omega)}^2+\lambda\|\phi(r)\|_{H^1(\bar{\Omega},d\sigma)}^2)\,dr+\|u(s)\|_{L^2(\Omega)}^2+\lambda\|\phi(s)\|_{H^1(\bar{\Omega},d\sigma)}^2.
\end{align}
From the classical Gronwall inequality, we deduce
\begin{align}\label{4.1.17}
\nonumber\|u(t+s)\|_{L^2(\Omega)}^2+\lambda\|\phi(t+s)\|_{H^1(\bar{\Omega},d\sigma)}^2\leq&(\|u(s)\|_{L^2(\Omega)}^2+\lambda\|\phi(s)\|_{H^1(\bar{\Omega},d\sigma)}^2)\exp(\int_s^{t+s}\mathbb{L}(r)\,dr)\\
\leq&\mathcal{M}_\ell(t)(\|u(s)\|_{L^2(\Omega)}^2+\lambda\|\phi(s)\|_{H^1(\bar{\Omega},d\sigma)}^2),
\end{align}
where
\begin{align}\label{4.1.18}
\mathcal{M}_\ell(t)=\exp(\int_0^{t+\ell}\mathbb{L}(r)\,dr)
\end{align}
is a finite number depending on $(u_{1_0},\phi_{1_0},\theta_{1_0})$ and $(u_{2_0},\phi_{2_0},\theta_{1_0})$ by using Theorem \ref{3.1}.

Integrating \eqref{4.1.17} with respect to $s$ for $0$ to $\ell,$ we obtain
\begin{align}\label{4.1.19}
\int_0^{\ell}\|u(t+s)\|_{L^2(\Omega)}^2+\lambda\|\phi(t+s)\|_{H^1(\bar{\Omega},d\sigma)}^2\,ds
\leq \mathcal{M}_\ell(t)\int_0^{\ell}\|u(s)\|_{L^2(\Omega)}^2+\lambda\|\phi(s)\|_{H^1(\bar{\Omega},d\sigma)}^2\,ds,
\end{align}
which implies the mapping $L_t:X_{\ell}\rightarrow X_{\ell}$ is locally Lipschitz continuous on $B_1^\ell$ for all $t\geq0.$\\
\hfill\qed

Thanks to the invariance of $B_1,$ Theorem \ref{4.1.1} and Theorem \ref{4.1.6}, we easily deduce that $K=\overline{L_{\tau_2} B_0^\ell}^{L^2(0,\ell;H\times V_I)}$ is positive invariant, uniformly absorbing compact subset of $X_\ell.$ Therefore, we can immediately obtain the existence of a global attractor in $X_\ell$ from Lemma \ref{2.7} stated as follows.
\begin{theorem}\label{4.1.20}
Assume that $h\in L^2(\Omega)$ and $(H_1)$-$(H_2)$ hold. Then the semigroup $\{L_t\}_{t\geq0}$ generated by problem \eqref{1.1}-\eqref{1.7} possesses a global attractor $\mathcal{A}_\ell$ in $X_\ell$ and $e_t(\mathcal{A}_\ell)$ is uniformly bounded in $H\times V_I$ with respect to $t\in [0,1],$ where
\begin{align*}
e_t(\mathcal{A}_\ell)=\{e_t(\chi):\chi\in\mathcal{A}_\ell\}
\end{align*}
for any $t\in [0,1].$
\end{theorem}

In what follows, we prove the smooth property of the semigroup $\{L_t\}_{t\geq0}$ to estimate the fractal dimension of the global attractor $\mathcal{A}_\ell.$
\begin{theorem}\label{4.1.21}
Assume that $h\in L^2(\Omega)$ and $(H_1)$-$(H_2)$ hold, let $\chi^1$ and $\chi^2$ be two short trajectories belonging to $\mathcal{A}_\ell.$ Then there exists a positive constant $\kappa$ independent of $t$ such that for arbitrary $t\geq\ell,$ we have
\begin{align*}
\|L_t\chi^1-L_t\chi^2\|_Y^2\leq\kappa\mathcal{M}_\ell(t)\int_0^{\ell}\|\chi^1(r)-\chi^2(r)\|_{H\times V_I}^2\,dr,
\end{align*}
where $\mathcal{M}_\ell(t)$ is given in \eqref{4.1.18}.
\end{theorem}
\textbf{Proof.} For any $\chi^1,$ $\chi^2\in \mathcal{A}_\ell,$ let $(u_1(t+\tau),\phi_1(t+\tau))=L_t\chi^1,$ $(u_2(t+\tau),\phi_2(t+\tau))=L_t\chi^2$ and let $u=u_1-u_2,$ $\phi=\phi_1-\phi_2.$ Since $e_t(\chi^1)$ and $e_t(\chi^2)$ is uniformly bounded in $H\times V_I$ with respect to $t\in [0,1]$ for any $\chi^1,$ $\chi^2\in \mathcal{A}_\ell,$ from the proof of Theorem \ref{3.1}, we obtain
\begin{align}\label{4.1.22}
\nonumber&\frac{d}{dt}(\|u(t)\|_{L^2(\Omega)}^2+\lambda\|\phi(t)\|_{H^1(\bar{\Omega},d\sigma)}^2)+\lambda\|\phi_t(t)\|_{L^2(\Gamma)}^2+\nu\|\nabla u\|_{L^2(\Omega)}^2+\lambda\gamma\|\nabla\Delta\phi\|_{L^2(\Omega)}^2\\
\leq&\mathbb{L}(t)(\|u(t)\|_{L^2(\Omega)}^2+\lambda\|\phi(t)\|_{H^1(\bar{\Omega},d\sigma)}^2),
\end{align}
where
\begin{align*}
\mathbb{L}(t)=C(1+\|\phi_1\|_{H^2(\bar{\Omega},d\sigma)}^2+\|\phi_2\|_{H^2(\bar{\Omega},d\sigma)}^2+\|\nabla u_2\|_{L^2(\Omega)}^2+\|\nabla\mu_2\|_{L^2(\Omega)}^2).
\end{align*}
For any $t\geq\ell,$ integrating \eqref{4.1.22} from $t-s$ to $t+\ell$ with $s\in[0,\frac{\ell}{2}],$ we conclude
\begin{align*}
&\|u(t+\ell)\|_{L^2(\Omega)}^2+\lambda\|\phi(t+\ell)\|_{H^1(\bar{\Omega},d\sigma)}^2+\int_{t-s}^{t+\ell}\lambda\|\phi_t(\zeta)\|_{L^2(\Gamma)}^2+\nu\|\nabla u(\zeta)\|_{L^2(\Omega)}^2+\lambda\gamma\|\nabla\Delta\phi(\zeta)\|_{L^2(\Omega)}^2\,d\zeta\\
\leq&\int_{t-s}^{t+\ell}\mathbb{L}(\zeta)(\|u(\zeta)\|_{L^2(\Omega)}^2+\lambda\|\phi(\zeta)\|_{H^1(\bar{\Omega},d\sigma)}^2)\,d\zeta+\|u(t-s)\|_{L^2(\Omega)}^2+\lambda\|\phi(t-s)\|_{H^1(\bar{\Omega},d\sigma)}^2.
\end{align*}
It follows from the classical Gronwall inequality that
\begin{align}\label{4.1.23}
\nonumber&\|u(t+\ell)\|_{L^2(\Omega)}^2+\lambda\|\phi(t+\ell)\|_{H^1(\bar{\Omega},d\sigma)}^2+\int_{t-s}^{t+\ell}\lambda\|\phi_t(\zeta)\|_{L^2(\Gamma)}^2+\nu\|\nabla u(\zeta)\|_{L^2(\Omega)}^2+\lambda\gamma\|\nabla\Delta\phi(\zeta)\|_{L^2(\Omega)}^2\,d\zeta\\
\leq&\exp(\int_{t-s}^{t+\ell}\mathbb{L}(\zeta)\,d\zeta)(\|u(t-s)\|_{L^2(\Omega)}^2+\lambda\|\phi(t-s)\|_{H^1(\bar{\Omega},d\sigma)}^2).
\end{align}
For any $t\geq\ell$ and any $s\in[0,\frac{\ell}{2}],$ integrating \eqref{4.1.22} from $s$ to $t-s,$ we obtain
\begin{align*}
&\|u(t-s)\|_{L^2(\Omega)}^2+\lambda\|\phi(t-s)\|_{H^1(\bar{\Omega},d\sigma)}^2\\
\leq&\int_s^{t-s}\mathbb{L}(r)(\|u(r)\|_{L^2(\Omega)}^2+\lambda\|\phi(r)\|_{H^1(\bar{\Omega},d\sigma)}^2)\,dr+(\|u(s)\|_{L^2(\Omega)}^2+\lambda\|\phi(s)\|_{H^1(\bar{\Omega},d\sigma)}^2).
\end{align*}
We deduce from the classical Gronwall inequality that
\begin{align}\label{4.1.24}
\nonumber\|u(t-s)\|_{L^2(\Omega)}^2+\lambda\|\phi(t-s)\|_{H^1(\bar{\Omega},d\sigma)}^2
\leq&(\|u(s)\|_{L^2(\Omega)}^2+\lambda\|\phi(s)\|_{H^1(\bar{\Omega},d\sigma)}^2)\exp(\int_s^{t-s}\mathbb{L}(r)\,dr)\\
\leq&(\|u(s)\|_{L^2(\Omega)}^2+\lambda\|\phi(s)\|_{H^1(\bar{\Omega},d\sigma)}^2)\exp(\int_0^{t-s}\mathbb{L}(r)\,dr).
\end{align}
Combining \eqref{4.1.23} with \eqref{4.1.24}, we obtain
\begin{align*}
\int_0^\ell\nu\|\nabla u(t+\zeta)\|_{L^2(\Omega)}^2+\lambda\gamma\|\nabla\Delta\phi(t+\zeta)\|_{L^2(\Omega)}^2\,d\zeta
\leq&\exp(\int_0^{t+\ell}\mathbb{L}(\zeta)\,d\zeta)(\|u(s)\|_{L^2(\Omega)}^2+\lambda\|\phi(s)\|_{H^1(\bar{\Omega},d\sigma)}^2)\\
=&\mathcal{M}_\ell(t)(\|u(s)\|_{L^2(\Omega)}^2+\lambda\|\phi(s)\|_{H^1(\bar{\Omega},d\sigma)}^2).
\end{align*}
Integrating the above inequality over $(0,\frac{\ell}{2})$ with respect to $s,$ we obtain
\begin{align*}
\int_0^\ell\nu\|\nabla u(t+\zeta)\|_{L^2(\Omega)}^2+\lambda\gamma\|\nabla\Delta\phi(t+\zeta)\|_{L^2(\Omega)}^2\,d\zeta
\leq\frac{2\mathcal{M}_\ell(t)}{\ell}\int_0^{\frac{\ell}{2}}\|u(s)\|_{L^2(\Omega)}^2+\lambda\|\phi(s)\|_{H^1(\bar{\Omega},d\sigma)}^2\,ds.
\end{align*}
Thanks to $\mathcal{M}_\ell(t)$ is bounded for any fixed $t\in [\ell,S],$ we obtain
\begin{align*}
\int_0^\ell\nu\|\nabla u(t+\zeta)\|_{L^2(\Omega)}^2+\lambda\gamma\|\nabla\Delta\phi(t+\zeta)\|_{L^2(\Omega)}^2\,d\zeta\leq \frac{2\mathcal{M}_\ell(t)}{\ell}\int_0^\ell\|u(s)\|_{L^2(\Omega)}^2+\lambda\|\phi(s)\|_{H^1(\bar{\Omega},d\sigma)}^2\,ds.
\end{align*}
It follows from the Sobolev trace Theorem and Lemma \ref{2.4} that
\begin{align}\label{4.1.25}
\int_0^\ell\nu\|\nabla u(t+\zeta)\|_{L^2(\Omega)}^2+\lambda\gamma\|\phi(t+\zeta)\|_{H^2(\bar{\Omega},d\sigma)}^2\,d\zeta\leq \kappa_1\mathcal{M}_\ell(t)\int_0^\ell\|u(s)\|_{L^2(\Omega)}^2+\lambda\|\phi(s)\|_{H^1(\bar{\Omega},d\sigma)}^2\,ds.
\end{align}
Thanks to
\begin{align}\label{4.1.26}
\nonumber\|u_t\|_{V^*}\leq&\nu\|\nabla u\|_{{L^2(\Omega)}}+C\|\nabla u\|_{{L^2(\Omega)}}\|\nabla u_2\|_{{L^2(\Omega)}}+C\|\nabla u_1\|_{{L^2(\Omega)}}\|\nabla u\|_{{L^2(\Omega)}}\\
&+C\|\phi_1\|_{H^1(\bar{\Omega},d\sigma)}\|\nabla\mu\|_{{L^2(\Omega)}}+C\|\phi\|_{H^1(\bar{\Omega},d\sigma)}\|\nabla\mu_2\|_{{L^2(\Omega)}}
\end{align}
and
\begin{align}\label{4.1.27}
\|\phi_t\|_{(H^1(\bar{\Omega},d\sigma))^*}\leq C\|\nabla u\|_{L^2(\Omega)}\|\phi_1\|_{H^1(\bar{\Omega},d\sigma)}+C\|\nabla u_2\|_{L^2(\Omega)}\|\phi\|_{H^1(\bar{\Omega},d\sigma)}+\gamma\|\nabla\mu\|_{L^2(\Omega)},
\end{align}
we infer from Theorem \ref{4.1.6}, \eqref{4.1.25}-\eqref{4.1.27} that
\begin{align}\label{4.1.28}
\left(\int_0^\ell\|u_t(t+r)\|_{V^*}+\|\phi_t(t+r)\|_{(H^1(\bar{\Omega},d\sigma))^*}\,dr\right)^2\leq \kappa_2\mathcal{M}_\ell(t)\int_0^\ell\|u(s)\|_{L^2(\Omega)}^2+\lambda\|\phi(s)\|_{H^1(\bar{\Omega},d\sigma)}^2\,ds.
\end{align}
The proof of Theorem \ref{4.1.21} is completed.\\
\qed\hfill

From Lemma \ref{2.8}, Theorem \ref{4.1.20} and Theorem \ref{4.1.21}, we immediately obtain the following result.
\begin{theorem}\label{4.1.29}
Assume that $h\in L^2(\Omega),$ $(H_1)$-$(H_2)$ hold. Then the fractal dimension of the global attractor $\mathcal{A}_\ell$ in $X_\ell$ of the semigroup $\{L_t\}_{t\geq0}$ generated by problem \eqref{1.1}-\eqref{1.7} established in Theorem \ref{4.1.20} is finite.
\end{theorem}
\subsection{The existence of a global attractor in $H\times V_I$}
In this subsection, we prove the existence of a finite dimensional global attractor in $H\times V_I$ of the semigroup generated by problem \eqref{1.1}-\eqref{1.7}.
\begin{theorem}\label{4.2.1}
Assume that $h\in L^2(\Omega)$ and $(H_1)$-$(H_2)$ hold. Then the mapping $e_1:\mathcal{A}_\ell\rightarrow \mathcal{A}=e_1(\mathcal{A}_\ell)$ is Lipschitz continuous. That is, for any two short trajectories $\chi^1,$ $\chi^2\in\mathcal{A}_\ell,$ there exists a positive constant $\theta$ dependent on $\ell$ such that
\begin{align*}
\|e_1(\chi^1)-e_1(\chi^2)\|_{H\times V_I}^2\leq\theta\int_0^{\ell}\|\chi^1(r)-\chi^2(r)\|_{H\times V_I}^2\,dr.
\end{align*}
\end{theorem}
\textbf{Proof.} For any $\chi^1,$ $\chi^2\in \mathcal{A}_\ell,$ let $(u_1(t+\tau),\phi_1(t+\tau))=L_t\chi^1,$ $(u_2(t+\tau),\phi_2(t+\tau))=L_t\chi^2$ and let $u=u_1-u_2,$ $\phi=\phi_1-\phi_2.$ Thanks to $e_0(\chi^1)$ and $e_0(\chi^2)$ is uniformly bounded in $H\times V_I$ for any $\chi^1,$ $\chi^2\in \mathcal{A}_\ell,$ from the proof of Theorem \ref{3.1}, we obtain
\begin{align}\label{4.2.2}
\nonumber&\frac{d}{dt}(\|u(t)\|_{L^2(\Omega)}^2+\lambda\|\phi(t)\|_{H^1(\bar{\Omega},d\sigma)}^2)+\lambda\|\phi_t(t)\|_{L^2(\Gamma)}^2+\nu\|\nabla u\|_{L^2(\Omega)}^2+\lambda\gamma\|\nabla\Delta\phi\|_{L^2(\Omega)}^2\\
\leq&\mathbb{L}(t)(\|u(t)\|_{L^2(\Omega)}^2+\lambda\|\phi(t)\|_{H^1(\bar{\Omega},d\sigma)}^2),
\end{align}
where
\begin{align*}
\mathbb{L}(t)=C(1+\|\phi_1\|_{H^2(\bar{\Omega},d\sigma)}^2+\|\phi_2\|_{H^2(\bar{\Omega},d\sigma)}^2+\|\nabla u_2\|_{L^2(\Omega)}^2+\|\nabla\mu_2\|_{L^2(\Omega)}^2).
\end{align*}
For $s\in (0,\ell),$ we infer from the classical Gronwall inequality and \eqref{4.2.2} that
\begin{align}\label{4.2.3}
\nonumber\|u(\ell)\|_{L^2(\Omega)}^2+\lambda\|\phi(\ell)\|_{H^1(\bar{\Omega},d\sigma)}^2\leq&(\|u(s)\|_{L^2(\Omega)}^2+\lambda\|\phi(s)\|_{H^1(\bar{\Omega},d\sigma)}^2)\exp(\int_s^\ell\mathbb{L}(r)\,dr)\\
\leq&(\|u(s)\|_{L^2(\Omega)}^2+\lambda\|\phi(s)\|_{H^1(\bar{\Omega},d\sigma)}^2)\exp(\int_0^\ell\mathbb{L}(r)\,dr).
\end{align}
Integrating \eqref{4.2.3} over $(0,\ell),$ we obtain
\begin{align*}
\|u(\ell)\|_{L^2(\Omega)}^2+\lambda\|\phi(\ell)\|_{H^1(\bar{\Omega},d\sigma)}^2
\leq\frac{1}{\ell}\exp(\int_0^\ell\mathbb{L}(r)\,dr)\int_0^\ell(\|u(s)\|_{L^2(\Omega)}^2+\lambda\|\phi(s)\|_{H^1(\bar{\Omega},d\sigma)}^2)\,ds.
\end{align*}
Thanks to \eqref{4.1.18}, we know that
\begin{align*}
\mathcal{M}_\ell=\exp(\int_0^\ell\mathbb{L}(r)\,dr)<+\infty,
\end{align*}
which implies that the mapping $e_1:\mathcal{A}_\ell\rightarrow \mathcal{A}$ is Lipschitz continuous.\\
\hfill\qed
\begin{theorem}\label{4.2.4}
Assume that $h\in L^2(\Omega)$ and $(H_1)$-$(H_2)$ hold. Then the semigroup $\{S_I(t)\}_{t\geq0}$ generated by problem \eqref{1.1}-\eqref{1.7} possesses a global attractor $\mathcal{A}=e_1(\mathcal{A}_\ell)$ in $H\times V_I.$ Furthermore, the fractal dimension of the global attractor $\mathcal{A}$ is finite.
\end{theorem}
\textbf{Proof.} From Lemma \ref{2.9}, Theorem \ref{4.1.29} and Theorem \ref{4.2.1}, we know that $\mathcal{A}$ is compact and the fractal dimension of $\mathcal{A}$ is finite. As a result of $L_t\mathcal{A}_\ell=\mathcal{A}_\ell,$ we have
\begin{align*}
S_I(t)\mathcal{A}=S_I(t)e_1(\mathcal{A}_\ell)=e_1(L_t\mathcal{A}_\ell)=e_1(\mathcal{A}_\ell)=\mathcal{A}
\end{align*}
for any $t\geq 0.$  From the definition of $B_1,$ we deduce that for any bounded subset of $H\times V_I,$ there exists some time $\bar{t}=\bar{t}(B)$ such that for any $t\geq \bar{t},$ we have
\begin{align*}
S_I(t)B\subset B_1.
\end{align*}
Therefore, we only need to prove that
\begin{align*}
\lim_{t\rightarrow+\infty}dist_{H\times V_I}(S_I(t)B_1,\mathcal{A})=0.
\end{align*}
Otherwise, there exist some positive constant $\epsilon_0,$ some sequence $\{(u_n,\phi_n)\}_{n=1}^{\infty}\subset B_1$ and some $\{t_n\}_{n=1}^{\infty}$ with $t_n\rightarrow+\infty$ as $n\rightarrow+\infty$ such that
\begin{align}\label{4.2.5}
dist_{H\times V_I}(S_I(t_n)(u_n,\phi_n),\mathcal{A})\geq\epsilon_0.
\end{align}
From the definition of $B_1,$ we deduce that there exists $\chi_n\in B_0^\ell$ such that
\begin{align*}
(u_n,\phi_n)=e_0(\chi_n).
\end{align*}
Since $\{\chi_n\}_{n=1}^{\infty}$ is bounded in $X_\ell$ and $\mathcal{A}_\ell$ is a global attractor in $X_\ell$ of the semigroup $\{L_t\}_{t\geq 0}$ generated by problem \eqref{1.1}-\eqref{1.7}, there exist a subsequence $\{\chi_{n_j}\}_{n=1}^{\infty}$ of $\{\chi_n\}_{n=1}^{\infty}$ and a subsequence $\{t_{n_j}\}_{n=1}^{\infty}$ of $\{t_n\}_{n=1}^{\infty}$ such that
\begin{align*}
L_{t_{n_j}-\ell}\chi_{n_j}\rightarrow\chi\in\mathcal{A}_\ell\,\,\,\textit{in}\,\,\, X_\ell\,\,\,\textit{for}\,\,\textit{as}\,\,\,j\rightarrow+\infty.
\end{align*}
Thanks to the continuity of $e_1,$ we have
\begin{align*}
S_I(t_{n_j})(u_{n_j},\phi_{n_j})=e_1(L_{t_{n_j}-\ell}\chi_{n_j})\rightarrow e_1(\chi)\in\mathcal{A}\,\,\,\textit{in}\,\,\, H\times V_I\,\,\,\textit{as}\,\,\,j\rightarrow+\infty,
\end{align*}
which contradicts with \eqref{4.2.5}.\\
\qed\hfill
\section*{Acknowledgement}
 This work was supported by
the National Science Foundation of China Grant (11401459).

\bibliographystyle{elsarticle-template-num}
\bibliography{BIB}
\end{document}